\documentclass[11pt,reqno,a4paper]{amsart} 
\usepackage[utf8]{inputenc}

\usepackage[margin=1.1in]{geometry}
\usepackage[parfill]{parskip}

\usepackage{amsmath,amsthm,amssymb}
\usepackage{mathtools}
\usepackage{esint}
\usepackage{xfrac}
\usepackage{comment}
\usepackage{enumerate}
\usepackage{enumitem}
\usepackage{cleveref}
\usepackage[final]{pdfpages}
\usepackage{bbm}
\usepackage{tikz}
\usetikzlibrary {lindenmayersystems}

\usepackage[normalem]{ulem} %%to make the strikethrough, removable at the end

\newtheorem{theorem}{Theorem}[section]
\newtheorem{lemma}[theorem]{Lemma}
\newtheorem{proposition}[theorem]{Proposition}
\newtheorem{corollary}[theorem]{Corollary}
\newtheorem{fact}{Fact}

\newtheorem*{theorem*}{Theorem}

\theoremstyle{definition}
\newtheorem{definition}[theorem]{Definition}
\newtheorem{case}{Case}

\theoremstyle{remark}
\newtheorem{remark}[theorem]{Remark}
\newtheorem{example}[theorem]{Example}

\numberwithin{equation}{section}

\renewcommand{\H}{\mathcal{H}}
\newcommand{\R}{\mathbb{R}}
\newcommand{\N}{\mathbb{N}}

\newcommand{\1}{\mathbbm{1}}

\newcommand{\loc}{\mathrm{loc}}
\newcommand{\eps}{\varepsilon}
\newcommand{\Cl}{\mathrm{Cl}}

\DeclareMathOperator{\divr}{div}
\DeclareMathOperator{\supp}{supp}
\DeclareMathOperator{\dist}{dist}

\newcommand{\jet}{\mathrm{jet}}
\DeclareMathOperator{\diam}{diam}

\title[Approximation of divergence-free vector fields]{Approximation of divergence-free vector fields vanishing on rough planar sets}

\author{Giacomo Del Nin}
\address[G.\ Del Nin]{Max Planck Institute for Mathematics in the Sciences, Inselstrasse 22, 04103 Leipzig, Germany}
\email{giacomo.delnin@mis.mpg.de}

\author{Bian Wu}
\address[B.\ Wu]{Max Planck Institute for Mathematics in the Sciences, Inselstrasse 22, 04103 Leipzig, Germany}
\email{bian.wu@mis.mpg.de}

\date{September 2024}

\begin{document}

\begin{abstract}
    Given any divergence-free vector field of Sobolev class $W^{m,p}_0(\Omega)$ in a bounded open subset $\Omega \subset \R^2$, we are interested in approximating it in the $W^{m,p}$ norm with divergence-free smooth vector fields compactly supported in $\Omega$. We show that this approximation property holds in the following cases: For $p>2$, this holds given that $\partial \Omega$ has zero Lebesgue measure (a weaker but more technical condition is sufficient); For $p\le 2$, this holds if $\Omega^c$ can be decomposed into finitely many disjoint closed sets, each of which is connected or $d$-Ahlfors regular for some $d\in[0,2)$. This has links to the uniqueness of weak solutions to the Stokes equation in $\Omega$. For H\"older spaces, we prove this approximation property in general bounded domains. \vspace{4pt}
    

    \vspace{4pt}
    \noindent\textsc{Keywords:} Sobolev spaces, Besov spaces, divergence-free, Ahlfors regular, Stokes equation, trace, Sard.
    %\vspace{5pt}  
\end{abstract}

\maketitle

%\tableofcontents

\section{Introduction}

\subsection{Background}

Consider the space $W^{m,p}_{0} (\Omega)$ of Sobolev functions supported in an open domain $\Omega \subset \R^n$, and the closely related space $\tilde W^{m,p}_{0} (\Omega)$ defined by
\begin{align*}
    W^{m,p}_{0}( \Omega ) :=& \,W^{m,p}-\text{closure of } C_c^{\infty}(\Omega), \\ 
    \tilde W^{m,p}_{0}( \Omega ) :=& \,\{ u \in W^{m,p}(\R^n) \mid D^ju=0 \text{ on } \Omega^c, \text{ for } 0 \leq |j| \leq m-1 \}. 
\end{align*}
Here, one has to interpret $u=0$ on $\Omega^c$ up to a set of $W^{m,p}$-capacity zero (see Section \ref{sec:capacity} for more details). The equivalence between $W^{m,p}_{0} (\Omega)$ and $\tilde W^{m,p}_{0} (\Omega)$ has attracted much attention since Sobolev's fundamental paper \cite{Sobolev} for smooth domains. This was later studied by Beurling \cite{Beurling}, Deny \cite{Deny1950}, Burenkov \cite{Burenkov}, Polking \cite{Polking}, and Hedberg \cite{Hedberg1973,Hedberg1978}. In 1980, a groundbreaking work of Hedberg \cite{Hedberg1981} proved the equivalence in \textit{general open subdomains} of $\R^n$ for $m \in \N^+$, $1 < p < \infty$ and $n \in \N^+$. For any function $f \in \tilde W^{m,p}_0(\Omega)$, Hedberg constructed a cutoff function $\omega \in C^\infty_c ( \Omega )$, which depends on $\Omega$ and on $f$, such that $\|f - \omega f\|_{W^{m,p}}$ is arbitrarily small (see Theorem \ref{thm:hedberg}). Later, Netrusov \cite{netrusov1992,Hedberg2007} proved the equivalence for general Besov spaces and Lizorkin-Triebel spaces.

\subsection*{A divergence-free vector-valued analog}
Since the Sobolev spaces of divergence-free vector fields arise in fluid PDEs, the vectorial analog of the above problem on divergence-free vector fields has also received considerable attention, i.e., the relation between $W^{m,p}_{0,\text{div}}( \Omega )$ and $\tilde W^{m,p}_{0,\text{div}}( \Omega )$, defined by
\begin{align}
    W^{m,p}_{0,\text{div}}( \Omega ) :=& \,W^{m,p}-\text{closure of } 
        C_c^{\infty}(\Omega, \R^n) \cap \{ \divr u = 0 \},  \label{e:intro:6} \\ 
    \tilde W^{m,p}_{0,\text{div}}( \Omega ) :=& \,\{ u \in W^{m,p}(\R^n, \R^n) \mid D^ju=0 \text{ on } \Omega^c, \text{ for } 0 \leq |j| \leq m-1, \,\divr u = 0 \}.  \label{e:intro:8}
\end{align}
Of particular importance is the case of $W^{1,2}$ divergence-free vector fields, as this has a connection to the uniqueness of solutions to the Stokes equation (see Section \ref{subsec:stokes_intro} for some more details about this).
Surprising examples given by Heywood \cite{Heywood1976} and Ladyzhenskaya \& Solonnikov \cite{Ladyzhenskaya1976} show that $W^{1,p}_{0,\text{div}}( \Omega )$ and $\tilde W^{1,p}_{0,\text{div}}( \Omega )$ are not identical for some unbounded locally smooth domain. The case of bounded domains was studied by Lions \cite{Lions1969}, Heywood \cite{Heywood1976}, Ladyzhenskaya \& Solonnikov \cite{Ladyzhenskaya1976}, Temam \cite{temam2024navier}. The most general results along this line are the equivalence for locally Lipschitz domains in dimensions $n \geq 2$, $m = 1$ and $p = 2$. Later, in the case $n \geq 2$, $m \geq 1$ and $p \in [1, \infty)$, Bogovski\u{\i} \cite{Bogovski1979,Bogovski1980} introduced a singular operator inverting the divergence operator, to prove the equivalence in domains with finitely many connected components, each of which is star-shaped. This singular operator is now known as Bogovski\u{\i} operator. More recently, Wang and Yang \cite{Wang2008} showed this coincidence for $n=2,3$, $k=1$ and $p=2$ in bounded domains with boundary satisfying a segment property. For general bounded domains, the equivalence between these two spaces remains an outstanding open problem (see \cite[III.7,Section~III.4]{galdi}).

In dimension $n=2$, divergence-free vector fields can be written as the rotated gradient of a scalar potential. Using this observation, for domains whose complement has finitely many connected components, \v{S}ver\'{a}k \cite{Sverak1990} pointed out that one can deduce the equivalence of $W^{1,p}_{0,\text{div}}( \Omega )$ and $\tilde W^{1,p}_{0,\text{div}}( \Omega )$ from Hedberg's result in the scalar case (see Subsection \ref{subsec:sverak}). However, nothing is known for generic domains with complement containing infinitely many connected components in dimension $n \geq 2$ (see also \cite[III.7,Section~III.4]{galdi} for more references about this problem). In this work, we partially fill this gap for rough planar domains, and we also prove the H\"{o}lder space counterpart \textit{without any assumption on bounded domain} $\Omega$.

\subsection{Main results}
\subsubsection*{The Sobolev case}
Our first result is for Sobolev spaces $W^{m,p}$, with $m \in \N^+$ and $p>2$.

%\giacomo{The point on $\Omega$ and $K,\tilde K$}
%\giacomo{Suggestion: $S(\Omega)$ and remove $K=\Omega^c$ (later we write $\Omega^c=K\cup\tilde K$.)}
\begin{theorem}[Divergence-free approximation in $W^{m,p}$, $p>2$]\label{thm:approximation_Sobolev_p>2}
    Suppose $m \in \N^+, p>2$ and a bounded open domain $\Omega \subset \R^2$ satisfies $|S( \Omega^c )|=0$, with
\begin{equation}
    \begin{split}
    S( \Omega^c) = \{ &x \in \Omega^c \mid x = \lim_{k \rightarrow \infty} x_k, \{ x_k \}_k \subset \Omega^c,\\ &\text{$x$ and 
        $x_k$ are in different connected components of $\Omega^c$ for each $k$} \}.
    \end{split}
\end{equation}
    Then for any $u\in W^{m,p}(\R^2;\R^2)$, satisfying $\divr u=0$ and $D^j u=0$ on $\Omega^c$ for any $0 \leq |j| \leq m-1$, there exists a sequence $\{u_k\}_k\subset C^\infty_c(\R^2;\R^2)$ satisfying $\divr u_k=0$ and $\supp u_k \subset \Omega$, such that $u_k\to u$ in $W^{m,p}$ as $k \rightarrow \infty$. 
\end{theorem}
In particular, as $S(\Omega^c)\subset\partial\Omega$, the class of bounded open domains $\Omega$ with $|S( \Omega^c )| = 0$ contains those whose boundary has Lebesgue measure zero. Here, $D^j u=0$ is well-defined, since by the Morrey-Sobolev embedding all derivatives $D^ju$, $|j|\le m-1$, can be assumed to be continuous and thus defined everywhere.

Our second result deals with the case $m \in \N^+$, $p \leq 2$ for bounded open sets whose complement is locally Ahlfors-regular. We recall that, for a fixed $d\in[0,n]$, a set $K\subset\R^n$ is called  $d$-regular if there exist $c,C>0$ such that
\[
    c r^d\le\H^d(K\cap B_r(x))\le C r^d\quad\text{$\forall x\in K, r\in(0,\diam K]$},
\]
where $\H^d$ stands for the $d$-dimensional Hausdorff measure.

\begin{theorem}[Divergence-free approximation in $W^{m,p}$, $p\le 2$]\label{thm:approximation_Sobolev_ple2}
    Let $m \in \N^+$ and let $\Omega \subset \R^2$ be a bounded open domain such that $\Omega^c=K\cup\tilde K$, where $K,\tilde K$ are disjoint closed sets, $\tilde K$ is the only unbounded component of $\Omega^c$, and $K$ is a compact, $d$-regular set for some $d\in[0,2)$. Then for any $u\in W^{m,p}(\R^2;\R^2)$, satisfying $\divr u=0$ and $D^j u=0$ on $\Omega^c$ for every $0 \leq |j| \leq m-1$, there exists a sequence $\{u_k\}_k\subset C^\infty_c(\R^2;\R^2)$ satisfying $\divr u_k=0$ and $\supp u_k \subset \Omega$, such that $u_k\to u$ in $W^{m,p}$ as $k \rightarrow \infty$. 
\end{theorem}

Here, since $u$ has only Sobolev regularity and $p\le 2$, the meaning of ``$D^j u=0$ on $K$'' must be understood up to $C_{m-|j|,p}$-capacity null sets (see Section \ref{sec:capacity}). As an example, if $m=1$ this means the following:
\begin{itemize}
    \item If $p>2$, then $u$ is continuous, and $u=0$ pointwise on $K$;
    \item If $2-d<p\le 2$, then $u$ admits a representative which is well-defined out of a $C_{1,p}$-null set, and $u$ is zero $C_{1,p}$-quasi everywhere on $K$;
    \item If $p\le 2-d$, then $C_{1,p}(K)=0$. In this case no information is imposed on the values of $u$ on $K$.
\end{itemize}
Similar considerations apply to higher-order derivatives if $m\ge 2$.

By a simple cutoff argument (see Lemma \ref{lemma:intersection_property}), we obtain that the same conclusion holds for more general domains $\Omega$.

\begin{corollary}   \label{cor:approximation_Sobolev_ple2}
The conclusion of \Cref{thm:approximation_Sobolev_ple2} holds true, if $\Omega$ is bounded and $\Omega^c$ admits a decomposition $\Omega^c = \bigcup_{0 \leq i \leq I} K_i$ for some $I \in \N^+$ and a collection of closed sets $\{K_i\}_{0 \leq i \leq I}$ satisfying the following properties:
\begin{enumerate}[label=(\roman*)]%[leftmargin=*,align=left]
    \item $K_0$ is the only unbounded connected component of $\Omega^c$;
    \item $K_i \cap K_\iota = \varnothing$ for any $i \neq \iota$;
    \item For any $1 \leq i \leq I$, $K_i$ is connected or $d_i$-regular for some $d_i \in [0,2)$.
\end{enumerate}
\end{corollary}

\begin{remark}[The role of Ahlfors-regularity]
Our proofs of Theorem \ref{thm:approximation_Sobolev_p>2} and Theorem \ref{thm:approximation_Sobolev_ple2} rely on certain trace theorems for the restriction of Sobolev functions to the set $K$. This is the reason why we restrict to sets satisfying $|S(\Omega^c)|=0$ in Theorem \ref{thm:approximation_Sobolev_p>2} and to Ahlfors-regular sets $K$ in Theorem \ref{thm:approximation_Sobolev_ple2}. Indeed, in these cases trace theorems are available. More precisely, when $K$ is $d$-regular the trace of a function in $W^{\alpha,p}(\R^n)$ on $K$ belongs to the Besov space $B^{p,p}_\beta(K)$, for $\beta=\alpha-\frac{n-d}{p}$. This is also the reason why we do not include the case $d=2$, as the corresponding trace theorem does not hold (see, however, Remark \ref{rmk:triebel_lizorkin}). We will recall the relevant theory from the monograph \cite{Jonsson-Wallin} in Section \ref{subsec:trace_besov}. We also provide a transparent, alternative proof of a trace theorem in the case $p>n$, which is a special case of a result in \cite{shvartsman2017whitney}. Proving trace theorems for more general sets $K$ would likely lead to a proof of the approximation property (namely, the analogue of \Cref{thm:approximation_Sobolev_p>2} and \Cref{cor:approximation_Sobolev_ple2}) for the corresponding sets. We also cite \cite{Vodopyanov-Tyulenev} for a trace theorem on $d$-thick sets (a more general class than $d$-regular sets), which however is not sufficient for our purposes, since we need a version for the space $W^{m,p}$, $m\ge 2$. 
\end{remark}

\subsubsection*{The H\"{o}lder case}

In addition to the theorems above, we also prove the corresponding theorem for vector fields with $C^{m,\gamma}$ regularity. In this case, we are able to treat any $m \ge 1$ and $\gamma \in [0,1]$ for \textit{all bounded open domains}. Up to our knowledge, this is the first result for proving the approximation property of divergence-free vector fields in general bounded domains \textit{without any regularity}. We also show that this result is sharp, providing in Example \ref{example:sharpness} a counterexample for the approximation property for vector fields with $C^{0,\gamma}$ regularity.

We present first the statement for $C^1$ vector fields separately, and we will give its proof in \Cref{ms:quicktour}. Besides having an interest in itself, this result also provides a ground to present, in a simplified way, the ideas that will be used for the proof of the Sobolev case. %With the aid of the added regularity, in $u$ is $C^2$ then we are able to prove the approximation theorem for any compact set $K$.

\begin{theorem}[Divergence-free approximation in $C^1$]\label{thm:approximation_C1}
    Let $K \subset \R^2$ be any compact set and $u \in C_c^2(\R^2;\R^2)$ with $\divr u = 0$. Suppose that $u=0$ and $\nabla u = 0$ on $K$. Then there exists a sequence $\{u_k\}_k\subset C^\infty_c(\R^2;\R^2)$ satisfying $\divr u_k=0$ and $\supp u_k\subset \R^2\setminus K$, such that $u_k \to u$ in $C^1$.
\end{theorem}

In order to state the result for general H\"{o}lder spaces we introduce the notation for the H\"{o}lder seminorm of a function $f$ on a set $E$, namely
%$|F|_{C^{0,\gamma}(K)}$ the H\"{o}lder norm of $F$ on $K$, namely
\[
    |f|_{C^{0,\gamma}(E)}:=\sup_{\substack{x,y\in E\\x\ne y}} \frac{|f(x)-f(y)|}{|x-y|^\gamma}.
\]
The full H\"{o}lder norm in an open set $U$ is then given by $\|f\|_{C^{m,\gamma}(U)}:=\|f\|_{C^m(U)}+|\nabla^m f|_{C^{0,\gamma}(U)}$.
Moreover, we denote by $K_\eps$ the $\eps$-neighborhood of $K$.

\begin{theorem}[Divergence-free approximation in $C^{m,\gamma}$] \label{thm:approximation_Cmgamma}
    Let $m \ge 1$ be an integer, and $0 \le \gamma \le 1$. Let $K \subset \R^2$ be any compact set and let $u \in C_c^{m,\gamma}(\R^2;\R^2)$ with $\divr u = 0$. 
    Suppose that
    \begin{equation}\label{eq:assumption_zero_derivatives_Cmgamma_main}
    D^j u(x)=0\quad\text{for every $x\in K$, $0\le |j|\le m$}.
    \end{equation}
    % \begin{equation}\label{eq:assumption_F_C^m-gamma}
    % \lim_{y\to x}\frac{|F(x)-F(y)|}{|x-y|^{m+\gamma}}=0 \quad\text{\giacomo{uniformly in $x\in K$}}
    % \end{equation}
    If $\gamma>0$ suppose in addition that
    \[
    \lim_{\eps\to 0^+}|\nabla^m u|_{C^{0,\gamma}(K_\eps)}=0.
    \]
    
    %\begin{equation}\label{eq:assumption_F_C^m-gamma_bis_main}
    %\frac{|\nabla^m u(x)-\nabla^m u(y)|}{|x-y|^\gamma}=0 \quad\text{\giacomo{uniformly in $x\in K$.}}
    %\end{equation}}

    % Suppose that
    % \[
    % \lim_{y\to x} \frac{|u(x)-u(y)|}{|x-y|^{m+\gamma}}=0\qquad\text{uniformly for $x\in K$},
    % \]
    % \giacomo{
    % and that
    % \[
    % \lim_{y\to x} \frac{|\nabla^m u(x)-\nabla^m u(y)|}{|x-y|^{\gamma}}=0\qquad\text{uniformly for $x\in K$}.
    % \]
    % }
    Then there exists a sequence $\{u_k\}_k\subset C^\infty_c(\R^2;\R^2)$ satisfying $\divr u_k=0$ and $\supp u_k\subset \R^2\setminus K$, such that $u_k \to u$ in $C^{m,\gamma}$.
\end{theorem}

\subsection{Connection to the Stokes operator}\label{subsec:stokes_intro}

The Stokes equation in a bounded domain $\Omega \subset \R^n$ is defined by
\begin{equation}   \label{e:stokes}
    \begin{split}
    - \Delta u + \nabla p =& f,  \quad \text{in } \Omega \\
    \divr u =& 0,   \quad \text{in } \Omega \\
    u =& 0,   \quad \text{on } \partial \Omega.
    \end{split}
\end{equation}
We say that $u \in W^{1,2}_0(\Omega, \R^n)$ is a weak solution of \eqref{e:stokes} if $\divr u = 0$ and for any $v \in C_c^\infty(\Omega, \R^n)$ with $\divr v = 0$, we have
\begin{align*}
    \int_\Omega \nabla u \cdot \nabla v \,dx
        = \int_\Omega f \cdot v \,dx.
\end{align*}
A crucial observation is the following: If $W^{1,2}_{0,\text{div}}( \Omega )$ and $\tilde W^{1,2}_{0,\text{div}}( \Omega )$ defined in \eqref{e:intro:6} and \eqref{e:intro:8} are not identical, then there exists a nonempty linear subspace $X \subset \tilde W^{1,2}_{0,\text{div}}( \Omega )$ such that 
\begin{align}
    \tilde W^{1,2}_{0,\text{div}}( \Omega )
        = X \oplus W^{1,2}_{0,\text{div}}( \Omega ).
\end{align}
Therefore, $X$ exactly contains all solutions to \eqref{e:stokes} with $f = 0$, which also means that the Stokes equation \eqref{e:stokes} \textit{always} has infinitely many solutions, even for $f \in C_c^\infty(\Omega, \R^n)$. From standard interior estimates for the Stokes operator (see \cite[Chapter IV]{galdi}), one can show that any function in $X$ belongs to $C^\infty(\Omega')$ for every open subdomain $\Omega'$ with $\bar \Omega' \subset \Omega$. Therefore, applying Riesz representation theorem, \Cref{thm:approximation_Sobolev_ple2} has the following corollary.

\begin{corollary}   \label{cor:stokes}
If $\Omega$ satisfies the condition specified in \Cref{cor:approximation_Sobolev_ple2}, then \eqref{e:stokes} admits a unique weak solution in $W^{1,2}_0(\Omega)$ with $f=0$.
\end{corollary}

\subsection{The main challenge and our approach}
The main difficulty is the nonlocality of the divergence-free condition, i.e. the product of a divergence-free vector field and a cutoff function may not be divergence-free. Therefore, contrary to the scalar case considered in Hedberg's works \cite{Hedberg1973,Hedberg1978,Hedberg1981}, the approach of multiplying with cutoff functions, which dates back to Beurling \cite{Beurling}, seems to be insufficient.

The works \cite{Bogovski1979,Bogovski1980} deal with divergence-free vector fields in Lipschitz or star-shaped domains using PDE methods, i.e. they involve solving the divergence equation in some form. This helps to correct the nonzero divergence introduced by using a cutoff function. However, the divergence equation is not always solvable in general open domains.

In this work we take a different approach, inspired by the observation made by \v{S}ver\'ak \cite{Sverak1990}. In dimension $2$, instead of looking at the vector field, we look at its scalar potential. In the case where $\Omega^c$ has finitely many connected components, it is straightforward to reduce the problem to Hedberg's result, as \v{S}ver\'ak noted in \cite{Sverak1990}. In this work, we introduce new ideas in dealing with general open domains. This helps us to resolve the approximation problem completely for H\"older based spaces in \Cref{thm:approximation_C1} and \Cref{thm:approximation_Cmgamma}. For Sobolev spaces, our approach needs a quantitative description of the trace of Sobolev functions onto compact domains. These types of results are only known in limited cases, which is the main obstruction to completely resolve the problem for Sobolev spaces.

\subsection{Outline of the paper} In \Cref{ms:pre} and \Cref{ms:hed} we recall or prove some preliminary results, including Whitney extension theorems, Morse-Sard theorems, Hedberg's theorem and a proof of \v{S}ver\'ak's observation. In \Cref{ms:quicktour} we prove the $C^1$ case \Cref{thm:approximation_C1} as \textit{a quick tour to our main ideas}. In \Cref{ms:trace} we prove \Cref{thm:approximation_Cmgamma}, the trace theorems and related results in Sobolev spaces. Finally, we present the proof of \Cref{thm:approximation_Sobolev_p>2}, \Cref{thm:approximation_Sobolev_ple2} and \Cref{cor:approximation_Sobolev_ple2} in \Cref{ms:main_proof}.

\subsection{Notation}%\giacomo{The notation section is almost useless now. I introduced the lebesgue measure below, so i think we can erase it.}
We will denote the Lebesgue measure of a set $A$ by $|A|$, and the $d$-dimensional Hausdorff measure of a set $E\subset \R^n$ by $\H^d(E)$. The ball of center $x$ and radius $r$ is denoted by $B(x,r)$ or $B_r(x)$. 
We use $\lesssim$ to denote the inequality $\leq$ up to a constant only depending on the dimension $n$, or on the fixed quantities within the proof, and similarly for $\gtrsim$. We use $\sim$ when both $\lesssim$ and $\gtrsim$ hold. $C_{s,p}(E)$ refers to the $(s,p)$-capacity of a set $E$ (see Section \ref{sec:capacity}). We use $\| \cdot \|_{p}$ to denote the $L^p$-norm. The open $\delta$-neighborhood of a set $A$ will be denoted by $A_\delta$, namely $A_\delta:=\{y:\, \dist(y,A)<\delta\}$. Given a real number $\beta$, we use $[\beta]$ to denote its integer part.

\subsection*{Acknowledgements}

The authors thank professors Hongjie Dong and Vladimir \v{S}ver\'ak for suggesting this problem during the special year program at Institute for Advanced Study from 2021 to 2022. Bian Wu thanks Institute for Advanced Study for the hospitality and the funding from Simons Foundation no. 816048. The authors also thank Hongjie Dong for locating the reference \cite{Weng} and Marjorie Drake for locating the reference \cite{shvartsman2017whitney}. The authors thank Hongjie Dong, Marjorie Drake, Max Goering, Lukas Koch and Stefan Schiffer for relevant discussions. 

\section{Preliminary tools and notation} \label{ms:pre}

In this section, we introduce some preliminary objects: the Whitney extension procedure, the notion of capacity, and the Morse-Sard theorem.
%\bian{We intruduce the following objects: Sobolev spaces on sets. Maximal functions. Whitney extension procedure.}

\subsection{Whitney decomposition, jets and extension} 
Referring to \cite[Chapter~VI]{Stein}, we recall some notation and results concerning Whitney-type extension theorems. We fix a compact set $K\subset \R^n$. Let $\{Q_k\}_k$ be a Whitney covering of $K^c$, namely a family of closed cubes with disjoint interiors, satisfying:
\begin{enumerate}[label=(\roman*)]
    \item $K^c = \bigcup_k Q_k$;
    \item $\sqrt{n}\ell(Q_k)\le \dist(Q_k,K)\le 4\sqrt{n}\ell(Q_k)$;
    \item If the boundaries of two cubes $Q_k$ and $Q_{k'}$ intersect, then $\frac14\le \frac{\ell (Q_{k})}{\ell (Q_{k'})}\leq 4$;
    \item For a given $Q_{k}$ there exist at most $12^{n}$ other cubes of the family that touch it.
\end{enumerate}
We also fix an associated partition of unity, namely a family of smooth functions $\{\varphi_k\}_k$, with $\supp \varphi_k$ contained in a small neighborhood of $Q_k$ (so that its support intersects only the cubes that touch $Q_k$), and reference points $\{y_k\}_k\subset K$ with the property that $\dist(Q_k,K)=\dist(Q_k,y_k)$. Here, we also define the index set of the cubes touching $Q_k$, i.e. its \textit{neighbors}
\begin{align}
    N(k) = \{ k' \mid k' \neq k, Q_{k'} \cap Q_k \neq \varnothing \}.
\end{align}

Next we introduce the space $J^m(K)$ of $m$-th order jets on $K$. Every element $\vec{f}$ of $J^m(K)$ is a collection $\vec{f}=\{f^{(j)}\}_{|j|\le m}$, where $f^{(j)}:K\to \R$ are functions and $j$ is a multi-index, namely $j\in \{1,\ldots,n\}^{m'}$, $0\le m'\le m$.

Given a jet $\vec{f}=\{f^{(j)}\}_{|j|\le m}$, we define the $m$-th order polynomial expansion of $\vec{f}$ centered at $y$ by
\begin{align} \label{f:taylor_poly_rev}
	P^{(m)}_y \vec{f} (x):=\sum_{|j|\le m} \frac{1}{j!} f^{(j)}(y) (x-y)^j .
\end{align}

For $\gamma\in(0,1]$ we define the space $C^{m,\gamma}_{\jet}(K)$ as the family of all jets $\vec{f}\in J^m(K)$ for which there exists a constant $M>0$ such that, defining
\[
    R_j \vec{f} (x,y) :=
    f^{(j)}(x) - \sum_{|j+l|\le m} \frac{1}{l!} f^{(j+l)}(y) (x-y)^l,
\]
we have
\begin{align*}
    |f^{(j)}(x)| \le M\quad\text{and}\quad |R_j \vec{f} (x,y)|\le M |x-y|^{m+\gamma-|j|}\qquad\text{for all $x,y\in K$, $|j|\le m$}.
\end{align*}
The $m$-th order jet norm of $\vec{f}$ on $K$ is given by the smallest constant $M$ satisfying the conditions above denoted by $\|\vec{f} \,\|_{C^{m,\gamma}_{\jet}(K)}$. Equivalently,
\[
    \|\vec{f}\|_{C^{m,\gamma}_{\jet}(K)} = 
        \sup_{\substack{|j|\le m\\ x,y\in K,\,x\ne y}} 
        \max\left\{|f^{(j)}(x)|, \, \frac{|R_j \vec{f} (x,y)|}{|x-y|^{m+\gamma-|j|}}\right\}.
\]
If $\gamma=0$ we define the space $C^{m,0}_{\jet}(K)$ (also simply denoted by $C^m_{\jet}(K)$) as the family of all jets $\vec{f}\in J^m(K)$ for which there exists a constant $M>0$ such that for every $|j|\le m$
\begin{align*}
    |f^{(j)}(x)| \le M\quad\text{and}\quad \frac{|R_j \vec{f} (x,y)|}{|x-y|^{m-|j|}}\to 0\qquad\text{uniformly in $x,y\in K$}.
\end{align*}
The norm $\|\vec{f}\|_{C^{m,0}(K)}$ is again defined as the smallest constant $M$ that satisfies all the inequalities above.
For a given $\gamma\in(0,1]$, we define the H\"{o}lder seminorm by
\begin{equation}\label{eq:def_holder_seminorm}
|f|_{C^{0,\gamma}(K)}:=\sup_{x,y\in K,\, x\ne y} \frac{|f(x)-f(y)|}{|x-y|^\gamma}.
\end{equation}
Finally, we define the pointwise Whitney extension of the jet $\vec{f}$ by
\begin{align} \label{f:Whitney_ext_pointwise_prelims}
    E^{(m)} \vec{f}(x) = 
        \begin{cases}
            f^{(0)}(x)     & x \in K, \\ 
            \sum_k \varphi_k(x) P^{(m)}_{y_k} \vec{f}(x),  & x \notin K.
        \end{cases}
\end{align}
We can now state the main extension theorems for $C^m$ and $C^{m,\gamma}$ spaces.

%\giacomo{Complete and check statemetns}

\begin{theorem}[Whitney extension in $C^m$] 
\label{thm:Whitney_C_m} 
Given $K \subset \R^n$ compact, $m \geq 1$ and $\vec{f} \in J^m(K)$, suppose that for every multi-index $j$ with $|j| \leq m$ it holds
\[
|f^{(j)}(x)|\le M,\quad\text{for every $x\in K$}
\]
and
\begin{align}
	\frac{|R_j \vec{f}(x,y)|}{|x-y|^{m-|j|}} \rightarrow 0, 
		\quad \text{ as } |x-y| \rightarrow 0, \quad \text{uniformly in } x,y\in K.
            \label{e:t:Whitney_C_m:6}
\end{align}
Then $F := E^{(m)} \vec{f} \in C^m(\R^n)$ and $D^j F = f^{(j)}$ on $K$ for any $|j| \leq m$. Furthermore, 
\[
\|F\|_{C^{m}}\lesssim \| \vec{f}\, \|_{C^{m}_{\jet}(K)}.
\]
%\giacomo{Is the above enough? Otherwise restore version commented below}
%we have, for any $j$ with $|j| \leq m$,
%\begin{align}
	%\|D^j F\|_{C^0} \lesssim \max\left\{\sup_{x\in K}|f^{(j)}(x),\sup_{x,y \in K, x \neq y} \frac{|R_j \vec{f}(x,y)|}{|x-y|^{m-|j|}}\right\}.
%\end{align}
\end{theorem}

\begin{theorem}[Whitney extension in $C^{m,\gamma}$] 
\label{thm:Whitney_C_m_gamma} 
Given $K \subset \R^n$ compact, $m \geq 1,\, \gamma \in (0,1]$ and $\vec{f} \in J^m(K)$, suppose that $\| \vec{f} \, \|_{C^{m,\gamma}_{\jet}(K)}<\infty$, i.e., there exists a constant $M$ such that for every multi-index $j$ with $|j| \leq m$ it holds
\[
|f^{(j)}(x)|\le M,\quad\text{for every $x\in K$.}
\]
and
\begin{align*}
	|R_j \vec{f}(x,y)| \leq M |x-y|^{m-|j| + \gamma},
		\quad \text{ for every } x,y \in K.
\end{align*}
Then $F := E^{(m)} \vec{f} \in C^{m,\gamma}(\R^n)$ and $D^j F = f^{(j)}$ on $K$. Furthermore, we have that
\[
\|F\|_{C^{m,\gamma}}\lesssim \| \vec{f} \, \|_{C^{m,\gamma}_{\jet}(K)}.
\]

%we have, for any $j$ with $|j| \leq m$,
%\begin{align}
%	\|D^j F\|_{C^{0,\gamma}} \lesssim C(|j|).
%\end{align}
\end{theorem}

\subsection{Capacity and Hausdorff measure}\label{sec:capacity}
We recall some basic notions regarding capacity and its relation with the Hausdorff measure (see \cite[Section~2.6]{ziemer}).

Given a set $E\subset \R^n$ and $s>0$, $p>1$, we define the $(s,p)$-capacity (also called $C_{s,p}$-capacity) of $E$ as
\[
C_{s,p}(E):=\inf\{\|f\|_{L^p(\R^n)}:\, f* g_s\ge 1\text{ on $E$} \},
\]
where $g_s(x)$ is the Bessel kernel, satisfying $\hat g_s(\xi)=(2\pi)^{-n/2}(1+|\xi|^2)^{-s/2}$.

We will use the following relation between capacity and Hausdorff measure.

\begin{theorem}[{\cite[Theorem~2.6.16]{ziemer}}]\label{thm:capacity_vs_hausdorff}
    Let $s>0$, $p>1$, $s p\le n$, and $E\subset \R^n$. The following implications hold:
    \begin{enumerate}[label=(\roman*)]
        \item If $\H^{n-s p}(E)<\infty$ then $C_{s,p}(E)=0$;
        \item If $C_{s,p}(E)=0$ then $\H^{n-s p+\eps}(E)=0$ for any $\eps>0$.
    \end{enumerate}
\end{theorem}

Putting together \cite[Theorem~3.1.4]{ziemer} and \cite[Theorem~3.3.3]{ziemer} we have the following.

\begin{theorem}[Approximate continuity out of capacity-null sets]\label{thm:app_continuity_out_of_capacity_null}
    Let $m \in \N^+$ satisfy $m p<n$, and let $u\in W^{m,p}(\R^n)$. Then there exists a set $E$ with $C_{m,p}(E)=0$ such that the limit
    \[
    \tilde u(x):=\lim_{r\to 0} \fint_{B_r(x)}u(y)\, dy
    \]
    exists for all $x\in E^c$. Moreover, for every $x\in E^c$
    \[
    \lim_{r\to 0} \fint_{B_r(x)}|\tilde u(x)-u(y)|^p\, dy=0.
    \]
\end{theorem}

\subsection{Morse-Sard theorem}
Recall that we denote by $|E|$ the Lebesgue measure of a set $E$. The classical Morse-Sard theorem for $C^2$ functions reads as follows.

\begin{theorem}[Morse-Sard in $C^2$]\label{thm:Sard}
    Let $\psi:\R^2\to \R$ be a $C^2$ function. Let $Z_\psi:=\{x\in \R^2:\, \nabla \psi(x)=0\}$ be the set of its singular points. Then $|\psi(Z_\psi)|=0$.
\end{theorem}

%Deleted: The cited result holds for maps in $BV_2$, but we just need the version for Sobolev spaces, so we will slightly modify the statements to adjust them to our needs.

We also recall the Sobolev version of the Morse-Sard theorem due to Bourgain, Korobkov and Kristensen \cite{BKK}. Following these authors, given $\psi\in W^{2,p}_\loc(\R^2)$ (which we can assume to be continuous by the Morrey-Sobolev embedding), we define
    \[
    Z_{0\psi}:=\Omega\cap \bigcap_{\eps>0} \Cl_M(\{x\in \Omega:\,|\nabla \psi(x)|\le \eps \})
    \]
    where, for a measurable set $A$, $\Cl_M(A)$ defined by
    \[
        \Cl_M(A):=\left\{x:\, \limsup_{r\to 0^+} \frac{|A\cap B_r(x)|}{|B_r(x)|}>0\right\}
    \]
    is the measure-theoretic closure of $A$. Observe that $Z_{0\psi}$ does not depend on the representative chosen for $\nabla \psi$. We further define
    \[
    Z_{1\psi}:=\{x\in\Omega:\,\text{$\psi$ is differentiable at $x$ and $\nabla\psi(x)=0$}\}.
    \]
    We set $Z_\psi:=Z_{0\psi}\cup Z_{1\psi}$. Now we cite the following results.
    \begin{theorem}[Morse-Sard in $W^{2,1}(\R^2)$ {\cite[Theorem~4.1]{BKK}}]\label{thm:Sard_Sobolev}\phantom{.}
        If $\psi\in W^{2,1}_\loc(\R^2)$, then we have $|\psi(Z_\psi)|=0$.
    \end{theorem}
    \begin{theorem}[Image of $\H^1$-null sets {\cite[Corollary~3.2]{BKK}}]\label{thm:image_null_sets}
        If $\psi\in W^{2,1}_\loc(\R^2)$, and $E\subset \R^2$ is a set with $\H^1(E)=0$, then $|\psi(E)|=0$.
    \end{theorem}

    Now, we prove the following result, needed in the proof of Theorem \ref{thm:approximation_Sobolev_ple2}.
    
    \begin{lemma}\label{lemma:sard_powered}
        Given $\psi\in W^{2,p}_\loc(\R^2)$ for some $p>1$ and $K\subset \R^2$ a compact set, assume $\nabla \psi=0$ $C_{1,p}$-quasi everywhere on $K$. Then we have $|\psi(K)|=0$. 
    \end{lemma}

    \begin{proof}
        By the assumption and the approximate continuity of $\nabla\psi$ out of $C_{1,p}$-null sets in Theorem \ref{thm:app_continuity_out_of_capacity_null}, there exists a set $G\subset K$ satisfying $C_{1,p}(K\setminus G)=0$ (and thus $\H^1(K\setminus G)=0$ by Theorem \ref{thm:capacity_vs_hausdorff}) and such that for every $x\in G$ it holds 
        \[
        \lim_{r\to 0^+}\fint_{B_r(x)}|\nabla \psi(y)|\, dy=0.
        \]
        We claim that $G\subset Z_{0\psi}$. Indeed, if $x\in G$ then 
        \[
        0 = \lim_{r\to 0^+}\fint_{B_r(x)}|\nabla \psi(y)|\, dy
            \ge \lim_{r\to 0^+} \frac{\eps}{|B_r(x)|}|\{y\in B_r(x): \, |\nabla\psi(y)|\ge \eps\}|
        \]
        which implies that $x\in \Cl_M(\{x\in \Omega:\,|\nabla \psi(x)|\le \eps \})$ for every $\eps>0$, hence the claim follows.

        We deduce that
        \[
        |\psi(K)|\le |\psi(G)|+|\psi(K\setminus G)|\le |\psi(Z_{0\psi})|+|\psi(K\setminus G)|=0,
        \]
        where for the last equality we use Theorem \ref{thm:Sard_Sobolev} and Theorem \ref{thm:image_null_sets}.
    \end{proof}

    \begin{corollary}\label{cor:sard}
        Let $\psi\in W^{2,p}(\R^2)$ for some $p>1$, and let $K\subset\R^2$ be a connected compact set. Assume that $\nabla\psi=0$ $C_{1,p}$-quasi everywhere on $K$. Then $\psi$ is constant on $K$.
    \end{corollary}

    \begin{proof}
        By continuity of $\psi$, the set $\psi(K)$ is a connected subset of $\R$, and thus an interval. On the other hand, by Lemma \ref{lemma:sard_powered} we know that $|\psi(K)|=0$. The only possibility is that $\psi(K)$ is a singleton, or equivalently, that $\psi$ is constant on $K$.
    \end{proof}

\subsection{\v{S}ver\'ak's result}\label{subsec:sverak}
We now present a generalization of \v{S}ver\'ak's observation in \cite[Lemma 5.2]{Sverak1990}. 

\begin{theorem} \label{thm:sverak}
    Let $\Omega\subset\R^2$ be a bounded domain such that $\Omega^c$ has finitely many connected components. Let $m\ge1$ and $1< p\le\infty$. Then $W^{m,p}_{0,\divr}(\Omega)=\tilde W^{m,p}_{0,\divr}(\Omega)$.
\end{theorem}

A proof for $m=1$, $p=2$ can also be found in the PhD thesis of Weng \cite{Weng}. We use the same idea, with the difference of employing the Sobolev version of Sard's theorem referenced above. 

\begin{proof}
    We write $\Omega^c=K_1\cup\ldots\cup K_N$, where each $K_i$ is a connected component of $\Omega^c$. Let $u\in W^{m,p}_{0,\divr}(\Omega)$, and let $\Psi$ be a potential for $u$, namely $\Psi\in W^{m+1,p}$ and $\nabla^\perp \Psi=u$. By Corollary \ref{cor:sard}, $\Psi=c_i$ on $K_i$ for some constants $c_i$. As $\Omega^c$ has only a finite number of connected components, we can find a smooth function $h\in C^\infty_c(\R^2)$ with the property that $h=c_i$ on a neighborhood of $K_i$ for every $i=1,\ldots,N$. Then the function $\Phi:=\Psi-h$ attains value zero on $\Omega^c$. By Hedberg's theorem (Theorem \ref{thm:hedberg}) we can find a sequence $\Phi_k\in C^\infty_c(\Omega)$ with $\Phi_k\to \Phi$ in $W^{m+1,p}$. It follows that the functions defined by $\Psi_k:=\Phi_k+h$ satisfy $\Psi_k\to \Psi$ in $W^{m+1,p}$, and moreover they attain value $c_i$ in a neighborhood of $K_i$. It follows that the vector fields defined by $u_k:=\nabla^\perp \Psi_k$ are zero on a neighborhood of $\Omega^c$ and converge to $u$ in $W^{m,p}$. Hence $u_k$ provide an approximating sequence for $u$, thus $u\in W^{m,p}_{0,\divr}(\Omega)$. This concludes the proof.
\end{proof}

\section{Hedberg's theorem} \label{ms:hed}

Our argument needs the following result by Hedberg \cite{Hedberg1981,Adams-Hedberg}, which constitutes a scalar version of the approximation theorem.

\begin{theorem}[{\cite[Theorem~9.1.3]{Adams-Hedberg}}]\label{thm:hedberg}
    Let $m$ be a positive integer, $1<p<\infty$, and $F\in W^{m,p}(\R^n)$. Let $\Omega\subset\R^n$ be an arbitrary open set. Then the following statements are equivalent:
    \begin{enumerate}[label=(\roman*)]
        \item $D^j F|_{\Omega^c}=0$ for all multi-indices $j$, $0\le|j|\le m-1$;
        \item $F\in W^{m,p}_0(\Omega)$;
        \item For any $\eps>0$ and any compact $E\subset\Omega$ there is a function $\eta\in C^\infty_0(\Omega)$ such that $\eta=1$ on $E$, $0\le\eta\le1$, and $\|F-\eta F\|_{W^{m,p}}<\eps$.
    \end{enumerate}

\end{theorem}

We refer the reader to \cite{Adams-Hedberg} for a proof of the above theorem. The same result for $C^{m,\gamma}$ spaces holds. We include its proof below, as we do not find it in the literature.  Recall the definition of H\"{o}lder seminorm in \eqref{eq:def_holder_seminorm}.

\begin{theorem}[Approximation in $C^{m,\gamma}$]\label{thm:Hedberg_Cmgamma}
    Let $m\in\mathbb{N}$, $0\le \gamma\le 1$. Let $K \subset \R^n$ be any compact set and $F \in C_c^{m,\gamma}(\R^n)$. 
    Suppose that
    \begin{equation}\label{eq:assumption_zero_derivatives_Cmgamma}
    D^j F(x)=0\quad\text{for every $x\in K$, $0\le |j|\le m$}.
    \end{equation}
    % \begin{equation}\label{eq:assumption_F_C^m-gamma}
    % \lim_{y\to x}\frac{|F(x)-F(y)|}{|x-y|^{m+\gamma}}=0 \quad\text{\giacomo{uniformly in $x\in K$}}
    % \end{equation}
    If $\gamma>0$ suppose also that
    \begin{equation}\label{eq:assumption_vanishing_holder_norm}
        \lim_{\eps\to 0^+} |\nabla^m F|_{C^{0,\gamma}(K_\eps)}=0.
    \end{equation}

    %\begin{equation}\label{eq:assumption_F_C^m-gamma_bis}
    %\frac{|\nabla^m F(x)-\nabla^m F(y)|}{|x-y|^\gamma}\le \omega(|x-y|)\qquad \forall x,y\in K_\eps
    %\end{equation}
    %for some modulus of continuity $\omega(\cdot)$, satisfying $\omega(s)\to 0$ as $s\to 0^+$.
    
    Then there exists a sequence $\{F_k\}_k\subset C^\infty_c(\R^n)$ satisfying $\supp F_k\subset \R^n\setminus K$, such that $F_k \to F$ in $C^{m,\gamma}$. Moreover, $F_k$ can be chosen of the form $F_k=F(1-\eta_k)$, for a sequence of smooth cutoff functions $\eta_k$.
\end{theorem}

\begin{remark}
    Observe that if $F\in C^m$ (and thus $\gamma=0$) then condition \eqref{eq:assumption_vanishing_holder_norm} automatically holds.
    Instead, if $\gamma>0$ then one can not just ask \eqref{eq:assumption_zero_derivatives_Cmgamma} to hold, as 
    one can consider $K=\{0\}$ and construct a function $F$ such that $|\nabla^m F|_{C^{0,\gamma}(B_\eps(0))}\ge 1$ for every $\eps$.
    %the following simple example shows: let $F(x)=|x|^{m+\gamma}$ and $K=\{0\}$. 
    Then every function with support outside $\{0\}$ has $C^{m,\gamma}$-distance at least 1 from $F$.
\end{remark}

\begin{proof}
    First we rewrite condition \eqref{eq:assumption_vanishing_holder_norm} in the following way: there exists a modulus of continuity $\omega(\cdot)$, i.e. $\omega(s)\to 0$ as $s\to 0^+$ and $\omega$ is increasing, such that
    \begin{equation}\label{eq:assumption_F_C^m-gamma_bis}
    |\nabla^m F(x)-\nabla^m F(y)|\le \omega(\eps) |x-y|^\gamma\qquad\text{for every $x,y\in K_\eps$.}
    \end{equation}
    For every $\eps>0$ we fix a smooth function $\rho_\eps$ with $\rho_\eps=1$ on $K_{\eps/4}$, $\rho_\eps=0$ on $K_{\eps}^c$, and
    \begin{equation}\label{eq:derivatives_of_rho_eps}
    |\nabla^k\rho_\eps|\le C(k)\eps^{-k}\qquad\text{for every $k\ge 0$}.
    \end{equation}
    For instance the choice $\rho_\eps=\chi_{K_{\eps/2}}*\zeta_{\eps/10}$ works, where $\zeta$ is a bump function supported on $B_{1}(0)$, $\int\zeta=1$ and $\zeta_s(y):=s^{-n}\zeta(y/s)$. We claim that the approximating family $F_\eps:=F(1-\rho_\eps)$ satisfies the conclusion of the theorem.
    
    Observe that $F_\eps$ belongs to $C^{m,\gamma}$ because $\rho_\eps$ is smooth away from $K$. Moreover, it is supported away from $K$ by construction. Therefore we just need to check the convergence, namely, that $\|F\rho_\eps\|_{C^{m,\gamma}}\to 0$ as $\eps\to 0$.
    
    First we check that for every multi-index $j$, with $|j|\le m$, it holds $\|D^j(F\rho_\eps)\|_{C^0(K_\eps)}\to 0$. By the Leibniz rule, we just need to check that $\|D^\theta F D^{j-\theta}\rho_\eps\|_{C^0(K_\eps)}\to 0$ for every multi-index $\theta\subseteq j$. We now show by (backwards)  induction that 
    \begin{equation}\label{eq:claim_D_theta_F}
    \|D^\theta F\|_{C^0(K_\eps)}\lesssim \omega(\eps) \eps^{m-|\theta|+\gamma}\quad\text{for every $|\theta|\le m$}.
    \end{equation}
    %and for some modulus of continuity $\tilde\omega(\cdot)$.
    If $|\theta|=m$ then given any $y\in K_\eps$ we can consider a point $x\in K$ of minimum distance from $y$. Thus $|x-y|\le \eps$, and from \eqref{eq:assumption_zero_derivatives_Cmgamma} and \eqref{eq:assumption_F_C^m-gamma_bis} we deduce that
    \[
    |D^\theta F(y)|\le \omega(\eps)\eps^\gamma.
    \]
    This proves the case $|\theta|=m$ of \eqref{eq:claim_D_theta_F}. Now consider $|\theta|<m$ and suppose by (backwards) induction that \eqref{eq:claim_D_theta_F} holds for every multi-index $\theta'$ with $m\ge|\theta'|>|\theta|$. Then for $x,y\in K_\eps$ such that $[x,y]\subset K_\eps$ it holds
    \begin{align}
    \begin{aligned}\label{eq:D_theta_bound}
        |D^\theta F(y)-D^\theta F(x)|&=\left|\int_0^1 \nabla (D^\theta F(x+t(y-x))\cdot (x-y)\, dt \right|\\
        &\lesssim \omega(\eps)\eps^{m-|\theta|-1+\gamma} |x-y|.
    \end{aligned}
    \end{align}
    In particular, by choosing again $x\in K$ as the point of minimum distance from $y$, and using that $D^\theta F$ vanishes on $K$, we obtain
    \begin{align*}
    |D^\theta F(y)|\lesssim \omega(\eps)\eps^{m-|\theta|+\gamma}.
    \end{align*}
    This proves the induction step and hence \eqref{eq:claim_D_theta_F} is verified.

    Recalling now also \eqref{eq:derivatives_of_rho_eps} it follows that
    \[
    \|D^\theta F D^{j-\theta}\rho_\eps\|_{C^0(K_\eps)}\lesssim \omega(\eps)\eps^{m-|\theta|+\gamma} \eps^{-|j|+|\theta|}=\omega(\eps)\eps^{m-|j|+\gamma},
    \]
    which is going to zero as $\eps\to 0$, since $|j|\le m$.

    We are now left to prove that $|\nabla^m(F\rho_\eps)|_{C^{0,\gamma}(K_\eps)}\to 0$, or equivalently, for every multi-index $j$ with $|j|=m$ and every multi-index $\theta\subseteq j$ we need to prove that
    \begin{equation}\label{eq:claim_holder_norm}
        |D^\theta FD^{j-\theta}\rho_\eps|_{C^{0,\gamma}(K_\eps)}\to 0.
    \end{equation}

    We claim the following bounds:
    \begin{enumerate}[label=(\roman*)]
        \item $\|D^\theta F\|_{C^0(K_\eps)}\le \omega(\eps)\eps^{m-|\theta|+\gamma}$; this follows directly from \eqref{eq:claim_D_theta_F}.
        \item $|D^\theta F|_{C^{0,\gamma}(K_\eps)}\lesssim \omega(2\eps)\eps^{m-|\theta|}$; to prove this, we fix $x,y\in K_\eps$ and consider two cases: if $|x-y|>\eps$ then
        \[
        \frac{|D^\theta F(y)-D^\theta F(x)|}{|x-y|^\gamma}\le \frac{2 \|D^\theta F\|_{C^0(K_\eps)}}{\eps^\gamma}\le 2\omega(\eps) \eps^{m-|\theta|}
        \]
        and we are done. If instead $|x-y|\le \eps$ then the segment $[x,y]$ is contained in $K_{2\eps}$, and reasoning as in \eqref{eq:D_theta_bound} we conclude that
        \[
        |D^\theta F(y)-D^\theta F(x)|\lesssim \omega(2\eps) \eps^{m-|\theta|-1+\gamma}|x-y|^\gamma\eps^{1-\gamma}=\omega(2\eps)\eps^{m-|\theta|}|x-y|^\gamma.
        \]
        This shows the claimed inequality.
        \item $\|D^{j-\theta}\rho_\eps\|_{C^0(K_\eps)}\lesssim \eps^{|\theta|-|j|}$; this follows directly from \eqref{eq:derivatives_of_rho_eps}.
        \item $|D^{j-\theta}\rho_\eps|_{C^{0,\gamma}(K_\eps)}\lesssim \eps^{|\theta|-|j|-\gamma}$; this follows from standard convolution estimates. Indeed,
        \begin{align*}
            |\rho_\eps(x)-\rho_\eps(y)|&\le \int_{B_\eps(x)\cup B_\eps(y)} \chi_{K_{\eps/2}}(z)|\zeta_{\eps/10}(x-z)-\zeta_{\eps/10}(y-z)|\,dz\\
            & \le \eps^{-n} \int_{B_\eps(x)\cup B_\eps(y)} \chi_{K_{\eps/2}}(z)|\zeta|_{C^{0,\gamma}} \left|\frac{x-y}{\eps}\right|^\gamma\,dz\\
            & \lesssim |x-y|^\gamma \eps^{-\gamma}.
        \end{align*}
        This shows the estimate when $|j|-|\theta|=0$, and for the general case one applies the same reasoning to the derivatives of $\rho_\eps$, keeping in mind that
        \[
        D^{j-\theta}\rho_\eps=\eps^{-|j|+|\theta|} (D^{j-\theta}\rho)_\eps.
        \]
    \end{enumerate}
    
    We are ready to prove \eqref{eq:claim_holder_norm}, also recalling the general bound 
    \[
    |fg|_{C^{0,\gamma}}\le \|f\|_{C^0}|g|_{C^{0,\gamma}}+|f|_{C^{0,\gamma}}\|g\|_{C^0}.
    \]    
    Using (i) and (iv) together, and then (ii) and (iii) together, we discover that both products are less than $\omega(2\eps)\eps^{m-|j|}$, which goes to zero as $\eps\to 0$.
    
    In conclusion, $\|F\rho_\eps\|_{C^{m,\gamma}(K_\eps)}\to 0$ and this finishes the proof.
    % \begin{align*}
    %     \frac{|F(x)\rho_\eps(x)-F(y)\rho_\eps(y)|}{|x-y|^{m+\gamma}}&= \frac{|F(x)(\rho_\eps(x)-\rho_\eps(y))|}{|x-y|^{m+\gamma}}+\frac{|(F(x)-F(y))\rho_\eps(y)|}{|x-y|^{m+\gamma}}\\
    %     &\le \|F\|_{C^0(K_{2\eps})}\frac{|\rho_\eps(x)-\rho_\eps(y)|}{|x-y|^{m+\gamma}}+\|\rho_\eps\|_{C^0}\frac{|F(x)-F(y)|}{|x-y|^{m+\gamma}}.
    % \end{align*}
    % The second term goes to zero by assumption as ...........................................\giacomo{finish proof}
\end{proof}

\begin{figure}
%\begin{minipage}{0.45\textwidth}
\begin{tikzpicture}
\draw [l-system={rule set={F -> F-F++F-F}, step=70pt, angle=-85,
   axiom=F, order=1}] lindenmayer system ;%-- cycle;
   \node at (0.5,2) {$\frac{\pi}{2}+\alpha$}; 
   \draw[->,gray] (0.8,1.7)--(2.2,0.3);
   \draw[->,gray] (0.9,1.8)--(3,0.3);
   
   \node at (4.5,1) {$2\alpha$}; 
   \draw[->,gray] (4,1)--(2.67,1.8);
   %\draw (2.67,2.45) circle[radius=0.05];
   %\draw (2.46,0) circle[radius=0.05];
   
   %%% ANGLES
   \draw[gray] ([shift=(-85:20pt)]2.67,2.45) arc (-85:-95:20pt);
   \draw[gray] ([shift=(85:10pt)]2.46,0) arc (85:180:10pt);
   \draw[gray] ([shift=(0:10pt)]2.88,0) arc (0:95:10pt);
   
   \draw[<->,gray] (0,-0.2)--(2.46,-0.2);
   \node at (1.25,-0.4) {$a$};

   \draw[<->,gray] (0,-0.6)--(5.35,-0.6);
   \node at (2.62,-0.8) {$1$};
\end{tikzpicture}\hspace{2cm}
%\end{minipage}
%\begin{minipage}{0.45\textwidth}
\begin{tikzpicture}
\draw [l-system={rule set={F -> F-F++F-F}, step=7pt, angle=-85,
   axiom=F, order=4}] lindenmayer system ;%-- cycle;
   \node at (2.62,-0.92) {};
\end{tikzpicture}
%\end{minipage}
\caption{The first and fourth step in the construction of the curve $g$. %The constant-speed parametrizations of these curves converge to $g$.
}\label{fig:koch}
\end{figure}
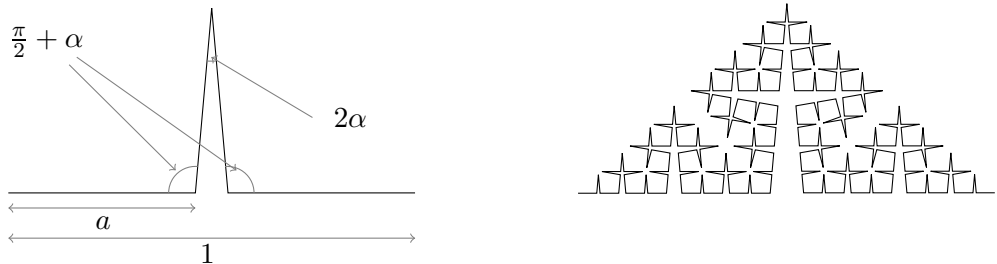

\begin{example}[Sharpness for H\"{o}lder spaces]\label{example:sharpness}
We shall now build an example which shows that for vector fields with $C^{0,\gamma}$ regularity, $\gamma\in[0,1)$, the approximation result analogous to Theorem \ref{thm:approximation_Cmgamma} is false. The idea is to use as potential $F$ a $C^{1,\gamma}$ function that is not constant on a connected set of critical points (the first such construction is due to Whitney \cite{whitney}), and then define the vector field by $\nabla^\perp F$. 

We take for granted the following fact, whose proof we could not find in the literature. We include a sketch of proof below, of which we do not claim any originality. 

    \begin{lemma}
        For every $\vartheta\in (\tfrac12,1]$ there exists a curve $g:[0,1]\to\R^2$ such that 
    \[
    C^{-1}|s-t|^\vartheta\le|g(s)-g(t)|\le C|s-t|^\vartheta
    \]
    for some $C>0$.
    \end{lemma}

    \begin{proof}[Sketch of proof]
        The construction mimics the iterative one used to build the Koch snowflake, but using a different angle between the segments (see Figure \ref{fig:koch}). We start from four consecutive segments of equal length $a>0$, forming consecutive exterior angles of size $\tfrac{\pi}{2}+\alpha$, $-2\alpha$, $\tfrac{\pi}{2}+\alpha$, for some $\alpha\in(0,\tfrac{\pi}{2}]$. We assume that the distance between the endpoints of the resulting curve is $1$. As a consequence, $a$ can attain any value in the interval $[\tfrac14,\tfrac12)$. At each step we replace each segment with a scaled-down copy of the previous step, the scaling factor being $a>0$. At step $k$ the resulting set identifies a non-self intersecting curve that can be parametrized with constant speed by some functions  $g_k:[0,1]\to\R^2$. One can show that $g_k$ form a Cauchy sequence in the $\vartheta$-H\"{o}lder norm, where $\vartheta=|\log_4 a|\in(\tfrac12,1]$. Thus they converge to some $g:[0,1]\to\R^2$. Moreover, one can show by induction that $|g_k(s)-g_k(t)|\ge C^{-1}|s-t|^\vartheta$ for some $C$ that does not depend on $k$. Hence this property passes to the limit $g$.
    \end{proof}

%that for every $\gamma\in[0,1)$ we can construct a $1/(1+\gamma)$-H\"{o}lder curve $g:[0,1]\to\R^2$ satisfying
%\[
%C^{-1}|s-t|^{1/(1+\gamma)}\le |g(t)-g(s)|\le C |s-t|^{1/(1+\gamma)}
%\]
%for some constant $C>0$.
Let us call $K:=g([0,1])$ and let us denote by $f:K\to [0,1]$ the inverse map of $g$, which thus satisfies
\[
|f(x)-f(y)|\le \tilde C |x-y|^{1+\gamma}\qquad\text{for every $x,y\in K$,}
\]
for some $\tilde C$, and where $\gamma:=\vartheta^{-1}-1$ can attain any value in $[0,1)$.

Let us consider the jet $\vec{f}=\{f^{(j)}\}_{|j|\le1}$ defined by:
\begin{align*}
f^{(0)}(x)&=f(x)\qquad\text{for every $x\in K$}\\
\qquad f^{(j)}(x)&=0\qquad\text{for every $x\in K$, $|j|=1$.}
\end{align*}
Then one can verify that $\vec{f}$ satisfies the conditions for Whitney's extension, namely Theorem \ref{thm:Whitney_C_m} or \ref{thm:Whitney_C_m_gamma}, and can thus be extended to a $C^{1,\gamma}$ function $F$ on $\R^2$. We can also assume that $F(0)=0$. Let us now consider the vector field $\bar u(x):=\nabla^\perp F(x)$. We claim that $\bar u$ can not be approximated in $C^{0,\gamma}$ by divergence-free vector fields supported out of $K$. 

To show this, assume by contradiction that there exists a sequence $u_k\in C^\infty_c(K^c)$ with $\divr u_k=0$ and $u_k\to u$ in $C^{0,\gamma}$. Denote by $F_k$ potentials for $u_k$, namely $C^{1,\gamma}$ functions such that $\nabla^\perp F_k=u_k$, $F_k(0)=0$. Then $F_k\to F$ in $C^0$: indeed for every $z\in \R^2$
\begin{align*}
|F_k(z)-F(z)|&=\left|F_k(0)+\int_{[0,z]}\nabla F_k(w)\cdot \tau(w)\,d\H^1(w)-F(0)-\int_{[0,z]}\nabla F(w)\cdot \tau(w)\,d\H^1(w)\right|\\
&\le\int_{[0,z]}|u_k(w)-u(w)|\,d\H^1(w)
\end{align*}
which is going to zero. However, $\nabla F_k=0$ on a neighborhood of $K$, and since $K$ is connected this implies that each $F_k$ is constant on $K$. This contradicts the fact that $F_k\to F$ in $C^0$, because $F$ is instead not constant on $K$.    
\end{example}

\section{A quick tour of the main idea: a proof in the $C^1$ case} \label{ms:quicktour}
In this section we give a full proof of Theorem \ref{thm:approximation_C1},  to showcase the main ideas for Sobolev spaces, namely Theorem \ref{thm:approximation_Sobolev_p>2} and Theorem \ref{thm:approximation_Sobolev_ple2}. We recall the statement for convenience:\\

%\giacomo{The statement for Holder spaces and Sobolev spaces is slightly different. }
%\begin{theorem} \label{C1_approximation}
\textit{Let $K \subset \R^2$ be a compact set and $u \in C_c^1(\R^2)$ with $\divr u = 0$. Suppose that $u=0$ and $\nabla u = 0$ on $K$. Then there exist $C^1$ vector fields $u_k$, compactly supported on $K^c$, such that $u_k \rightarrow u$ in $C^1$.}\\
%\end{theorem}

\begin{proof}[Proof of \Cref{thm:approximation_C1}]

We conduct this proof in several steps.

\begin{enumerate}[label=(S\arabic*)]
    \item \textit{Potential.} Since $\divr u=0$, we consider a potential $\Psi$ for $u$, namely, $\Psi \in C_c^2(\R^2)$ such that $\nabla^\perp \Psi=u$. First, we define
    \[
    \Psi^*(x_1,x_2):= \int_{0}^{x_1} -u_2(t,0) dt + \int_{0}^{x_2} u_1(x_1,s)ds.
    \]
    Then we define $\Psi = \Psi^* - c_0$ for some constant $c_0 \in \R$, such that $\Psi$ has compact support. In particular, $\nabla\Psi=0$ and $\nabla^2\Psi=0$ on $K$ and $u = \nabla^\perp \Psi$.

    \item \label{item:res} \textit{Restriction.} We can consider the restriction of $\Psi$ to $K$ as a jet of order 2:
    \begin{fact}\label{fact:restriction}
        $\Psi$ defines a family of second-order jets at each point by $\vec{\psi} = \{\psi^{(j)}\}_{|j|\le 2} :=\{D^j\Psi\}_{|j|\le 2}\in J^2(K)$. In particular, $\psi^{(0)} = \Psi|_K$ and  $\psi^{(j)}=0$ for $|j|\ge 1$.
    \end{fact}
    
    \item \label{item:morse} \textit{Morse-Sard.} Theorem \ref{thm:Sard} allows us to conclude the following fact.
    \begin{fact}\label{fact:sard}
        $\Psi(K)$ is a compact set with zero measure.
    \end{fact}
    For any $\eps>0$, we cover $\Psi(K)$ with finitely many open intervals $I_1^\eps,\ldots,I_N^\eps$ with disjoint closures and total measure less than $\eps$. Taking the preimages of $I_i^\eps$, we find disjoint open sets $U_1^\eps,\ldots,U_N^\eps\subset\Omega$ that are well-separated from one another, and such that $K\subset U_1^\eps\cup\ldots\cup U_N^\eps$. By shrinking $U_i^\eps$ for each $i$ if necessary, we can also assume that $U^\eps:=\bigcup_i U_i^\eps$ is contained in the $\eps$-neighborhood of $K$.

    \item \label{item:compression} \textit{Monotone compression.} We define a compression map $\eta_\eps:\R\to\R$ by
    \begin{align}   \label{e:eta_eps_qt}
    \eta_\eps(t):=\int_0^t \1_{\bigcup_i I_i^\eps }(s)\, ds.
    \end{align}
    This map has the following properties: 
    \begin{enumerate}
        \item There are real numbers $c_i^\eps$ such that $\eta_\eps(t)=t-c_i^\eps$ on $I_i^\eps$;
        \item $\|\eta_\eps\|_{C^0}\le \eps$;
        \item $ | \eta_\eps(t) - \eta_\eps(t) | \leq |t-s| $ for any $t,s \in \R$.
    \end{enumerate}
    Now we define the compressed function $\psi_\eps^{(0)}:=\eta_\eps\circ \Psi|_K$. This function satisfies $\psi_\eps^{(0)} = \psi^{(0)} - c_i^\eps$ on $K \cap U_i^\eps$, and $\|\psi_\eps^{(0)}\|_{C^0(K)} \le \eps$. We also set $\psi_\eps^{(j)}:=\psi^{(j)}=0$ for $1\le|j|\le 2$. We denote $\vec{\psi}_\eps = \{ \psi^{(j)}_\eps \}_{|j|\le 2} \in J^2(K)$.
    
    \item\label{item:compression_estimate} \textit{Estimate on the compressed jets.} We claim the following fact.
    \begin{fact}\label{fact:norm_goes_to_0}
        $\|\vec{\psi}_\eps\|_{C^2_{\jet}(K)}\to 0$ as $\eps\to 0$.
    \end{fact}
    Observe that
    \[
    \|\vec{\psi}_\eps\|_{C^2_{\jet}(K)} = \sup_{x,y\in K,\, x\ne y}
        \frac{|\psi_\eps^{(0)}(x)-\psi_\eps^{(0)}(y)|}{|x-y|^2}.
    \]
    To prove Fact \ref{fact:norm_goes_to_0}, assume by contradiction that it does not hold. Then there exists $\kappa > 0$, and sequences $\{(x_k,y_k)\}_k$ and $\{\eps_k\}_k$ such that
    \begin{align}
      x_k, y_k \in K,& \quad \eps_k \rightarrow 0 \\ 
      \frac {|\Psi(x_k)-\Psi(y_k)|} {|x_k-y_k|^{2}} &\geq 
            \frac {|\psi_{\eps_k}^{(0)}(x_k)-\psi_{\eps_k}^{(0)}(y_k)|} {|x_k-y_k|^{2}} 
                \geq \kappa.  \label{e:p:approximation_C1:12}
    \end{align}
    Here we used property \ref{item:compression}(c) for $\eta_{\eps_k}$. It follows from \eqref{e:p:approximation_C1:12} and \ref{item:compression}(b) that
    \begin{align*}   
      |x_k-y_k|^{2} \leq \frac{\eps_k} {\kappa} \rightarrow 0, \text{ as } k \rightarrow \infty.
    \end{align*}
    This contradicts $\nabla^2 \Psi = 0$ and $\nabla \Psi=0$ on $K$. Indeed, denoting by $[x_k,y_k]$ the segment between $x_k$ and $y_k$, we have
    \begin{align}   \label{e:p:approximation_C1:16}
    |\Psi(x_k)-\Psi(y_k)|\le \|\nabla^2\Psi\|_{C^0([x_k,y_k])} |x_k-y_k|^2,
    \end{align}
    and the segment $[x_k,y_k]$ is contained in the $\eps_k/\kappa$-neighborhood of $K$, thus 
    \begin{align}
     \|\nabla^2\Psi\|_{C^0([x_k,y_k])}
        \le \|\nabla u^\perp\|_{C^0(K_{\eps_k/\kappa})}\to 0. \label{e:p:approximation_C1:18}
    \end{align}
    Hence, with $\psi^{(0)} = \Psi|_K$, \eqref{e:p:approximation_C1:12}, \eqref{e:p:approximation_C1:16} and \eqref{e:p:approximation_C1:18} give a contradiction.

    \item\label{item:extension} \textit{Extension.} By Whitney's extension theorem we have the following.
    \begin{fact}\label{fact:extension}
        We can extend $\vec{\psi}_\eps$ to $\varphi_\eps \in C_c^2(\R^2)$ satisfying $\|\varphi_\eps\|_{C^2(\R^2)}
            \lesssim \|\vec{\psi}_\eps\|_{C^2_{\jet}(K)} = o(\eps)$. 
    \end{fact}
    First, we need to check the condition \eqref{thm:Whitney_C_m}. This follows from $\nabla^2 \Psi = 0$ and $\nabla \Psi=0$ on $K$ and the properties of $\eta_{\eps_k}$. Then, we apply \Cref{thm:Whitney_C_m} to get $E^{(2)} \vec{\psi}_\eps \in C^2(\R^2)$ with $\|E^{(2)} \vec{\psi}_\eps\|_{C^2(\R^2)}
            \lesssim \|\vec{\psi}_\eps\|_{C^2_{\jet}(K)}$.
    Since $K$ is compact and $U^\eps$ is bounded, one can multiply $E^{(2)} \vec{\psi}_\eps$ with a cutoff function to get $\varphi_\eps$ with compact support and $\|\varphi_\eps\|_{C^2(\R^2)}
            \lesssim \|\vec{\psi}_\eps\|_{C^2_{\jet}(K)}$. Therefore, we have \Cref{fact:extension}.
    
    From \ref{item:morse} and \ref{item:compression}, we also know that the function $\Psi-\varphi_\eps$ is constant on each $K\cap U_i^\eps$, with value $c_i^\eps$.

    \item\label{item:auxiliary} \textit{Auxiliary function.} Note that $\{K\cap U_i^\eps\}_i$ are finitely many disjoint closed sets, hence there exists small $\rho > 0$ and a function $h_\eps\in C_c^\infty(\R^2)$ such that $h_\eps=c_i^\eps$ on $K_\rho \cap U_i^\eps$. Now $g_\eps := \Psi-\varphi_\eps-h_\eps \in C_c^2(\R^2)$ with
    \begin{align*}
        \supp \nabla h_\eps \subset&\, K^c, \\
        g_\eps = 0,\, \nabla g_\eps =&\, 0, \,\nabla^2 g_\eps = 0,
            \quad \text{on } K.
    \end{align*}
    %Here $K_\delta$ denotes the $\delta$-neighbourhood of $K$.

    \item\label{item:Hedberg} \textit{Hedberg's theorem and conclusion.} We invoke Theorem \ref{thm:Hedberg_Cmgamma} to find a sequence $g_\eps^k\in C^\infty_c(\R^2\setminus K)$ with $\|g_\eps^k-g_\eps\|_{C^2}\to 0$ as $k\to\infty$, and we define the sequence $u_\eps^k:=\nabla^\perp(g_\eps^k+h_\eps)$. Since from \ref{item:auxiliary} we have that $\nabla^\perp(g_\eps^k+h_\eps)=0$ on some neighborhood of $K$, we deduce that $u_\eps^k \in C^\infty_c( K^c )$ . Also, we can estimate
    \begin{align*}
        \|u_\eps^k-u\|_{C^1} &= \|\nabla^\perp(g_\eps^k+h_\eps-\Psi)\|_{C^1} \\
        &\le \| \nabla^\perp (g_\eps+h_\eps-\Psi) \|_{C^1} + \|\nabla^\perp(g_\eps-g_\eps^k)\|_{C^1}\\
        &= \|\nabla^\perp \varphi_\eps\|_{C^1} + \| \nabla^\perp(g_\eps^k-g_\eps) \|_{C^1}\\
        &\le \|\varphi_\eps\|_{C^2}+\|g_\eps^k-g_\eps\|_{C^2}.
    \end{align*}
        Now the first summand goes to zero as $\eps\to 0$ by Fact \ref{fact:norm_goes_to_0}, while the second summand goes to zero as $k\to\infty$ by Theorem \ref{thm:Hedberg_Cmgamma}. Considering a diagonal sequence we can extract an approximating sequence $\{ u_k \}_k$ converging to $u$ in $C^1$, and this concludes the proof.\qedhere
    \end{enumerate}
\end{proof}

\begin{remark}
%\giacomo{Check if this remark is true or not...}
    % In the $C^2$ case we could actually avoid the construction of an auxiliary function and the use of Hedberg's theorem, and we could merge steps \ref{item:extension}),\ref{item:auxiliary}),\ref{item:Hedberg}) by means of a single extension. Indeed, it is not hard to show (with essentially the same proof) that Fact 2 holds in the stronger version
    % \[
    % \|\psi_\eps\|_{C^2_{\jet}(K_{\delta(\eps)})}\to 0\quad\text{as $\eps\to 0$},
    % \]
    % for some neighbourhoods $K_{\delta(\eps)}$ of $K$. By Whitney's extension we could thus extend $\psi|_{K_{\delta(\eps)}}$ to a $C^2$ function $\varphi_\eps$, and the sought approximating sequence would be $u_n:=\nabla^\perp(\psi-\varphi_\eps)$, since $\psi-\varphi_\eps$ is locally constant on $K_{\delta(\eps)}$. 
    We chose to write down the proof in the $C^1$ case separately, because the proof strategy remains almost the same for Sobolev spaces. To prove the approximation result for Sobolev spaces we shall replace Facts \ref{fact:restriction}-\ref{fact:extension} with the corresponding statements for Sobolev spaces. More precisely:
    \begin{itemize}
        \item Fact 1 is replaced by trace theorems for Sobolev functions which shall be addressed in \Cref{ms:trace}.
        \item Fact 2 is replaced by the corresponding Morse-Sard result for Sobolev functions proven by Bourgain, Korobkov and Kristensen \cite{BKK} given in Lemma \ref{lemma:sard_powered}.
        \item Fact 3 requires an argument that hinges upon the characterization of traces, with the $C^2_{\jet}$ norm replaced by the appropriate Besov norm, see Proposition \ref{prop:compression_besov} and Proposition \ref{prop:compression_sobolev}.
        \item Fact 4 is replaced by a corresponding extension theorem from Besov to Sobolev functions, see Theorem \ref{thm:extension_sobolev_p>2} and Theorem \ref{thm:trace_besov}.
    \end{itemize}
    
\end{remark}

%Deleted: The extension and trace theorems for Sobolev functions are the reason why we restrict to $d$-regular sets $K$ when $p\le n$, since this is the regularity required by the best Sobolev trace/extension results known to date \cite{Jonsson-Wallin}. We stress that the extension of Theorem \ref{thm:approximation_Sobolev_ple2} to an arbitrary compact set $K$ would likely follow if one had suitable trace/extension theorems for arbitrary compact sets. 

\section{Trace theorems and monotone compression}  \label{ms:trace}

In this section we prove or quote trace results for Sobolev spaces.
%(namely, Sobolev spaces with $p>n$, Sobolev spaces with $p\le n$ and H\"{o}lder spaces). 
We also present the argument of \textit{monotone compression} for H\"older spaces and Sobolev spaces, analogous to \ref{item:compression} and \ref{item:compression_estimate} in Section \ref{ms:quicktour}. For this purpose, we recall the compression maps $\eta_\eps$: given a compact set $K$ and a continuous function $\Psi$ with $|\Psi(K)|=0$, and given $\eps>0$, we consider a finite family of disjoint open intervals $I_1^\eps,\ldots ,I_N^\eps$ that cover $\Psi(K)$ and that satisfy $|I_1^\eps|+\ldots+|I_N^\eps|<\eps$. Then we set
\begin{equation}\label{eq:def_eta_eps}
    \eta_\eps(t):=\int_0^t \1_{\bigcup_i I_i^\eps }(s)\, ds.
\end{equation}
We also recall the following three properties from \ref{item:compression} in Section \ref{ms:quicktour}:
\begin{align}
    |\eta_\eps(t)-\eta_\eps(s)| \le&\, |t-s|,   \label{eq:eta_eps_1}   \\
    \|\eta_\eps\|_{C^0} \le&\, \eps,   \label{eq:eta_eps_2}\\
    \text{There are real numbers $c_i^\eps$}&\text{ such that $\eta_\eps(t)=t-c_i^\eps$ on $I_i^\eps$} \label{eq:eta_eps_3}.
\end{align}

%\subsection{Trace theorem in H\"{o}lder spaces}

%For H\"{o}lder spaces, we fix the pointwise Whitney extension of the jet $\vec{f}$ to be
%\begin{align} \label{f:Whitney_ext_pointwise}
%    E^{(m)} \vec{f}(x) = 
%        \begin{cases}
%            f^{(0)}(x)     & x \in K, \\ 
%            \sum_k \varphi_k(x) P^{(m)}_{y_k} \vec{f}(x),  & x \notin K.
%        \end{cases}
%\end{align}

\begin{lemma}[Monotone compression in H\"{o}lder spaces]\label{lemma:compression_Holder}
    Let $F\in C^{m,\gamma}(\R^n)$ and $K$ a compact set satisfying
    \begin{equation}\label{eq:assumption_zero_derivatives_Cmgamma_body}
    D^j F(x)=0\quad\text{for every $x\in K$, $1\le |j|\le m$}.
    \end{equation}
    % \begin{equation}\label{eq:assumption_F_C^m-gamma}
    % \lim_{y\to x}\frac{|F(x)-F(y)|}{|x-y|^{m+\gamma}}=0 \quad\text{\giacomo{uniformly in $x\in K$}}
    % \end{equation}
    If $\gamma>0$ suppose in addition that
    \begin{equation}\label{eq:assumption_holder_compression_m}
    \lim_{\eps\to 0^+}|\nabla^m F|_{C^{0,\gamma}(K_\eps)}=0.
    \end{equation}
    Let $\vec{f}=\{f^{(j)}\}_{|j|\le m}$ be the $m$-th order jet given by $f^{(j)}:=D^j F|_K$. Let $\eta_\eps$ be the compression map defined in \eqref{eq:def_eta_eps}. Then the family $\vec{f}_\eps=\{ f_\eps^{(j)}\}_{|j|\le m}$ defined by 
\[
 f_\eps^{(j)}:=
\begin{cases}
    \eta_\eps\circ f^{(0)} &\text{if $j=0$}\\
    0 & \text{otherwise}
\end{cases}
\]
satisfies
\begin{equation}\label{eq:smallness_compression_holder}
    \|\vec{f}_\eps\|_{C^{m,\gamma}_{\jet}(K)}\le \delta(\eps)
\end{equation}
for some function $\delta:\R_+\to\R_+$ with $\delta(\eps)\to 0$ as $\eps\to 0^+$.
\end{lemma}

\begin{proof}
    Observe that by assumption $f^{(j)}=0$ for $|j|\ge 1$, hence we have 
    \begin{equation} \label{eq:holder_ratio}
        \| \vec f_\eps \|_{C^{m,\gamma}_{\jet}(K)} := 
        \max \bigg\{ \| f_\eps^{(0)} \|_{C^0(K)}, \,
        \sup_{x,y\in K,\, x\ne y} \frac{ |f_\eps^{(0)}(x)-f_\eps^{(0)}(y)| }{|x-y|^{m+\gamma}} 
        \bigg\}.
    \end{equation}
    We assume by contradiction that the conclusion of the lemma does not hold. Since we know $\|f_\eps\|_{C^0(K)} \to 0$ as $\eps \to 0$, the last term in \eqref{eq:holder_ratio} does not converge to zero as $\eps\to 0$. This entails the existence of $\kappa > 0$, and of sequences $\{(x_k,y_k)\}_k$ and $\{\eps_k\}_k$ such that
    \begin{align}
      x_k, y_k \in K,& \quad  \eps_k \rightarrow 0 \\ 
      \frac {|F(x_k)-F(y_k)|} {|x_k-y_k|^{m+\gamma}} 
            &\geq \frac {|f^{(0)}_{\eps_k}(x_k)-f^{(0)}_{\eps_k}(y_k)|} {|x_k-y_k|^{m+\gamma}} \geq \kappa.
    \end{align}
    It follows from the properties of $\eta_{\eps_k}$ that
    \[
      |x_k-y_k|^{m+\gamma} \leq \frac{\eps_k} {\kappa} \rightarrow 0, \text{ as } k \rightarrow \infty.
    \]
    This contradicts \eqref{eq:assumption_holder_compression_m} and \eqref{eq:holder_ratio}. Indeed, denoting by $[x_k,y_k]$ the segment between $x_k$ and $y_k$, we can estimate
    \[
        |F(x_k) - F(y_k)| \le \|\nabla F\|_{C^0([x_k,y_k])} |x_k-y_k|.
    \]
    Moreover for every $z\in [x_k,y_k]$, if $|j|< m$, then we have
    \[
    |D^j F(z)|=|D^j F(x_k)-D^j F(z)|\le C \|\nabla^{|j|+1}F\|_{C^0([x_k,y_k])} |x_k-y_k|.
    \]
    If $|j|=m$, then
    \[
    |D^j F(z)|=|D^j F(x_k)-D^j F(z)|\le |\nabla^m F|_{C^{0,\gamma}([x_k,y_k])} |x_k-y_k|^\gamma.
    \]
    Putting together the estimates above, we obtain for $\gamma=0$ that
    \[
        |f(x_k)-f(y_k)| \lesssim |x_k-y_k|^{m}\|\nabla^mF\|_{C^0([x_k,y_k])},
    \]
    and for $\gamma>0$ that
    \[
        |f(x_k)-f(y_k)| \lesssim |x_k-y_k|^{m+\gamma}|\nabla^m F|_{C^{0,\gamma}([x_k,y_k])}.
    \]
    Since the segment $[x_k,y_k]$ is contained in the $\eps_k/\kappa$-neighborhood of $K$, by either \eqref{eq:assumption_zero_derivatives_Cmgamma_body} if $\gamma=0$ or \eqref{eq:assumption_holder_compression_m} if $\gamma>0$, we have that
    \[
    \frac{|f(x_k)-f(y_k)|}{|x_k-y_k|^{m+\gamma}}
        \rightarrow 0, \quad \text{as $k \to \infty$.}
    \]
    This gives a contradiction, which concludes the proof.
\end{proof}

Now we can prove \Cref{thm:approximation_Cmgamma}.

\begin{proof}[Proof of \Cref{thm:approximation_Cmgamma}]
    The proof is the same as the proof of \Cref{thm:approximation_C1}, apart from the following technical points. If $u\in C^{m,\gamma}_c(\R^2;\R^2)$, with $m\ge 1$, $\gamma\in[0,1)$, then the potential $\Psi$ belongs to $C^{m+1,\gamma}(\R^2)$. We replace Fact \ref{fact:restriction} and Fact \ref{fact:extension} with Whitney's extension given by Theorem \ref{thm:Whitney_C_m} and Theorem \ref{thm:Whitney_C_m_gamma}. By Lemma \ref{lemma:compression_Holder} applied with $F=\Psi$,  $\|\vec{\psi}_\eps\|_{C^{m,\gamma}_{\jet}(K)}\to 0$ as $\eps\to 0$. Finally, in Point \ref{item:Hedberg} one uses Theorem \ref{thm:Hedberg_Cmgamma} with a similar reasoning.
\end{proof}

\subsection{Trace theorem in Sobolev spaces with $p>n$}

In this subsection we give a transparent proof of a special case of the trace theorems by Shvartsman \cite{shvartsman2017whitney}. A full trace description has two components, the extension part and the restriction part. We recall the following maximal function for jets defined by Shvartsman,
\begin{align} \label{maximal_Shvartsman}
      M^{(m)} \vec{f} (x) := \sup_{y,z \in K, y \neq z} \frac{ |P^{(m-1)}_{y} \vec{f} (x) - P^{(m-1)}_{z} \vec{f} (x)| } { |x-y|^m + |x-z|^m },
\end{align}
and the definition of the Hardy–Littlewood maximal function $M f$ for $f \in L^1_{\text{loc}}(\R^n)$\
\begin{align*}
    Mf(x) = \sup_{r>0} \frac{1}{r^n} \int_{B_r(x)} |f(y)| dy.
\end{align*}

% \begin{theorem}[Trace for $W^{m,p}$, $p>n$]
%     Let $K \subset \R^n$ be compact, $m \geq 1,\, p > n$ and let $\vec{f} := \{f^{(j)}\}_{|j| \leq m-1} \in J^{m-1}(K)$ with $f^{(j)} = 0$ for any $|j| \ge 1$. Then $\vec{f}$ is the trace of A $W^{m,p}$ function if and only if $M^{(m)}\vec{f}\in L^p(\R^n)$. More precisely:
%     \begin{itemize}
%         \item \textbf{Restriction:} Let $F\in W^{m,p}(\R^n)$. Then for every multi-index $j$ with $|j|\le m-1$ the functions $D^jF$ admit a continuous representative and thus $f^{(j)}:=D^jF|_K$ are pointwise defined.  and satisfy
%         \item \textbf{Extension:} 
%     \end{itemize}
% \end{theorem}

\begin{theorem}[Extension in $W^{m,p}$ for $p>n$]
\label{thm:extension_sobolev_p>2}
Given $K \subset \R^n$ compact, $m \geq 1,\, p > n$ and $\vec{f} := \{f^{(j)}\}_{|j| \leq m-1} \in J^{m-1}(K)$ with $f^{(j)} = 0$ for every $|j| \ge1$, suppose
\begin{align*}
	\| M^{(m)}\vec{f} \, \|_{L^p} < +\infty.
\end{align*}
Then $F := E^{(m-1)} \vec{f} $ belongs to $W^{m,p}(\R^n)$, and for every $j$ with $1\le |j| \leq m-1$, we have $D^j F|_K = 0$. Furthermore, we have 
\begin{align}
	\|\nabla^m F\|_{L^p} \lesssim \| M^{(m)}\vec{f} \, \|_{L^p}.
\end{align}

\end{theorem}

\begin{proof}
Define
\begin{align} \label{f:Whitney_ext_pointwise_der}
    H^{(j)} \vec{f}(x) := 
        \begin{cases}
            f^{(j)}(x)     & x \in K, \\ 
            \displaystyle\sum_{j=l+\theta} D^l \varphi_k(x) D^\theta P^{(m-1)}_{y_k} \vec{f}(x),  
                & x \notin K.
        \end{cases}
\end{align}
In this proof, we only work with $(m-1)$-jets, so we introduce a short-hand notation $P_z = P^{(m-1)}_z$.
We would like to prove that, for any $l$ with $|l| \leq m-1$,
\begin{align}   \label{e:Whitney_trace_Sobolev_holder:7}
    \| \nabla H^{(l)} \vec{f} \|_{L^p} \lesssim \| M^{(m)}\vec{f} \, \|_{L^p},
\end{align}
and $D^j H^{(l)} \vec{f} = H^{(j+l)} \vec{f}$ when $|j|=1$. It suffices to prove, for fixed $1 \leq i \leq n$, any $h \in \R$ with $0 < h < \varepsilon_U$, any $j,l$ with $|j|=1$ and $|l|, |l+j| \leq m-1$, we have
\begin{align}
    \| \tau_h H^{(l)} \vec{f} - H^{(l)} \vec{f} \, \|_{L^p} \lesssim \, h \| M^{(|l|+1)}\vec{f} \,& \|_{L^p},
        \quad \tau_h F(x) := F(x+h e_i),    \label{e:Whitney_trace_Sobolev_holder:8} \\
    \Big| \frac{ \tau_h H^{(l)} \vec{f}(x) - H^{(l)}\vec{f}(x) } {h} - H^{(j+l)} \vec{f}(x) \, \Big|
        \lesssim& \, h M^{(|l|+2)}\vec{f}(x) , 
        \quad \text{for any } x \in K.  \label{e:Whitney_trace_Sobolev_holder:9}
\end{align}
Here, $\varepsilon_U$ is a small universal constant. Indeed, \eqref{e:Whitney_trace_Sobolev_holder:8} implies that $H^{(l)} \vec{f} \in W^{1,p}(\R^n)$. Note that $H^{(l)} \vec{f}$ is smooth in $K^c$, then \eqref{e:Whitney_trace_Sobolev_holder:9} implies the difference quotient at $x$ in \eqref{e:Whitney_trace_Sobolev_holder:9} converges pointwise to $H^{(j+l)}\vec{f}(x)$ for every $x \in \R^n$. By dominated convergence, we have that $D^j H^{(l)} \vec{f} = H^{(j+l)} \vec{f}$ belongs to $L^p$, and this concludes the proof. It thus remains to show \eqref{e:Whitney_trace_Sobolev_holder:8} and \eqref{e:Whitney_trace_Sobolev_holder:9}. 

We divide the proof in four different cases.

\begin{case}[$x, x + he_i \in K$]
Taking $y=x$ and $z=x + he_i$. This is obvious from the definitions of the maximal function $M^{(\cdot)}$ for jets.
\end{case}

\begin{case}[$x+he_i \in K, x \in K^c$ or $x \in K, x+he_i \in K^c$] We consider the case $x+he_i \in K, x \in K^c$. The other one is analogous. Suppose $x \in Q_{k_0}$. Hence, for any $k \in N(k_0)$, $|x-y_k| \lesssim h$. %We also choose $\varepsilon_U$ small such that $|x+he_i-y_k| \leq \frac{1}{4}$ for any $k \in N(k_0)$. 
When $|l|=0$,
\begin{align*}
	\frac{|\tau_h F(x) - F(x)|} {h} 
        =&\, \frac{1}{h} \Big| F(x+he_i) - \sum_{k \in N(k_0)} \varphi_k(x) f^{(0)}(y_k) \Big| \\ 
		\lesssim&\, \frac{1}{h} \sum_{k \in N(k_0)} \varphi_k(x) \big| f^{(0)}(x+he_i) - f^{(0)}(y_k) \big| \\ 
        \lesssim& \, h \min \big\{ M^{(2)}(x), M^{(2)}(x+he_i) \big\}.
\end{align*}
For $|l|\geq 1$, we have
\begin{align*}
    \frac{|\tau_h H^{(l)} \vec{f}(x) - H^{(l)} \vec{f}(x)|} {h} 
        =&\, \frac{1}{h} \Big| \sum_{k \in N(k_0)} D^l \varphi_k(x) P_{y_k} \vec{f}(x) \Big| \\ 
        \lesssim&\, \frac{1}{h} \Big| \sum_{k \in N(k_0)} D^l \varphi_k(x) \Big( f^{(0)}(y_k) - f^{(0)}(y_{k_0}) \Big) \Big|\\
        \lesssim&\, h \min \big\{ M^{(|l|+2)}(x), M^{(|l|+2)}(x+he_i) \big\}.
\end{align*}
The above implies \eqref{e:Whitney_trace_Sobolev_holder:8} and \eqref{e:Whitney_trace_Sobolev_holder:9} with $H^{(j+l)} \vec{f} = 0$.
\end{case}

\begin{case}[$x \in Q_{k_0}$, $x+h \in Q_{k_1}$, and $|h| \geq \frac{1}{10} \min \{ \ell(Q_{k_0}), \ell(Q_{k_1}) \}$]
As above, for any $k \in N(k_0)$ and any $k' \in N(k_1)$, $|x-y_k| + |x-y_{k'}| \lesssim |h|$.
Then for $|l|=0$, we have
\begin{align*}
	\frac{|\tau_h F(x) - F(x)|} {h} =&\, \frac{1}{h} \Big| \sum_{k \in N(k_0)} \varphi_k(x) P_{y_k} \vec{f}(x) - \sum_{k' \in N(k_1)} \varphi_{k'}(x) P_{y_{k'}} \vec{f}(x+h) \Big| \\ 
		\leq&\, \frac{1}{h} \sup_{k \in N(k_0), k' \in N(k_1)} 
		\big| f^{(0)}(y_{k}) - f^{(0)}(y_{k'}) \big| 
        \lesssim\, h M^{(2)}(x).
\end{align*}
For $|l|\geq 1$, we have
\begin{align*}
    \frac{|\tau_h H^{(l)} \vec{f}(x) - H^{(l)} \vec{f}(x)|} {h} 
        =&\, \frac{1}{h} \Big| \sum_{k \in N(k_0)} D^l \varphi_k(x) P_{y_k} \vec{f}(x) - \sum_{k' \in N(k_1)} D^l \varphi_{k'}(x+h) P_{y_{k'}} \vec{f}(x+h) \Big| \\ 
        \lesssim& \frac{1}{h} \sup_{k \in N(k_0), k' \in N(k_1)} \ell(Q_{k_0})^{-|l|} \big| f^{(0)}(y_k) - f^{(0)}(y_{k'}) \big| \\ 
        \lesssim& \, h \min \big\{ M^{(|l|+2)}(x), M^{(|l|+2)}(x+he_i) \big\}.
\end{align*}
These give \eqref{e:Whitney_trace_Sobolev_holder:8} and \eqref{e:Whitney_trace_Sobolev_holder:9} with $H^{(j+l)} \vec{f} = 0$.

\end{case}

\begin{case}[$x \in Q_{k_0}$, $x+h \in Q_{k_1}$, and $|h| \leq \frac{1}{10} \min \{ \ell(Q_{k_0}), \ell(Q_{k_1}) \}$] For any $t \in [0,1]$, $x+th \in K^c$, and we compute
\begin{align*}
    |\nabla H^{(l)} \vec{f}(x+th)| =&\, \Big| \sum_{k \in N(k_0) \cup N(k_1)} \nabla D^l \varphi_k(x+th) f^{(0)}(y_k) \\ 
        &- \sum_{k \in N(k_0) \cup N(k_1)} \nabla D^l \varphi_k(x+th) f^{(0)}(y_{k_0}) \Big| \\ 
        \lesssim& \sum_{k \in N(k_0) \cup N(k_1)} |\nabla D^l \varphi_k(x+th)| | f^{(0)}(y_k) - f^{(0)}(y_{k_0}) | \\
        \lesssim& \sum_{k \in N(k_0) \cup N(k_1)} \ell(Q_{k_0})^{-(|l|+1)} | f^{(0)}(y_k) - f^{(0)}(y_{k_0}) | \\ 
        \lesssim&\, M^{(|l|+1)}(x) \lesssim h M^{(|l|+2)}(x).
\end{align*}
Notice that $\tau_h H^{(l)} \vec{f}(x) - H^{(l)} \vec{f}(x) = \int_0^1 h \cdot \nabla H^{(l)} \vec{f}(x+th) dt$. Then 
\begin{align*}
    \frac{|\tau_h H^{(l)} \vec{f}(x) - H^{(l)} \vec{f}(x)|} {|h|} \lesssim M^{(|l|+1)}(x) \lesssim h M^{(|l|+2)}(x).
\end{align*}
\end{case}

\end{proof}

\begin{theorem}[Restriction in $W^{m,p}$ for $p>n$]\label{thm:restriction_sobolev_p>2}
Given $K \subset \R^n$ compact, for $F \in W^{m,p}(\R^n)$, define a jet $\vec{f} \in J^{m-1}(K)$ by $f^{(j)} = D^j F$ on $K$ for any $|j| \leq m-1$. Then
\begin{align*}
    \| M^{(m)}\vec{f} \, \|_{L^p} \lesssim \|\nabla^m F\|_{L^p}.
\end{align*}
\end{theorem}

\begin{proof}
We recall the following inequality from Mazya \cite{mazya2011}. For fixed $q \in (n,p)$, for any $|j| \leq m-1$, any cube $Q$ and any $x,y \in Q$, we have
\begin{equation} \label{Sobolev_inequality_pointwise}
    \big| D^j F(x) - D^j P_y^{(m-1)} \vec{f} (x) \big|
            \lesssim \ell(Q)^{m-|j|} \left( \fint_{Q} |\nabla^m F|^q \right)^{1/q}.
\end{equation}
Then for any $Q$ containing $y,z \in K$, we have
\begin{align*}
    \big| D^j P_y^{(m-1)} \vec{f} (x) - D^j P_z^{(m-1)} \vec{f} (x) \big|
        \lesssim \ell(Q)^{m-|j|} \left( \fint_{Q} |\nabla^m F|^q \right)^{1/q}.
\end{align*}
Indeed, we can obtain the inequality above by applying \eqref{Sobolev_inequality_pointwise} twice for $x,y$ and $x,z$ respectively. Now let $Q$ be centered at $x$ with radius $|x-y|+|x-z|$ and $j=0$. We obtain
\begin{align*}
    \frac{\big| P_y^{(m-1)} \vec{f} (x) - P_z^{(m-1)} \vec{f} (x) \big|} {|x-y|^m+|x-z|^m}
        \lesssim&\, \left( \fint_{Q} |\nabla^m F|^q \right)^{1/q} \\
        \lesssim&\, \Big( M \big( |\nabla^m F|^q \big)(x) \Big)^{1/q}.
\end{align*}
Then we have
\begin{align*}
    \| M^{(m)}\vec{f} \, \|_{L^p}^p 
        \lesssim&\, \int \Big( M \big( |\nabla^m F|^q \big)(x) \Big)^{p/q} dx \\ 
        \lesssim&\, \big\| M \big( |\nabla^m F|^q \big) \big\|_{L^{p/q}}^{p/q} 
        \lesssim\, \| \nabla^m F \|_{L^p}.\qedhere
\end{align*}
\end{proof}

Next, we show that the monotone compression \ref{item:compression} in \Cref{ms:quicktour} works well in the Sobolev space $W^{m,p}(\R^2)$.

\begin{proposition}[Monotone compression for Sobolev spaces, $p>2$] \label{prop:compression_sobolev}
Let $m \geq 2$, $p>1$, and consider a compact set $K \subset \R^2$ that satisfies $|S(K)|=0$, with
\begin{equation}
    \begin{split}
    S(K) = \{ &x \in K \mid x = \lim_{k \rightarrow \infty} x_k, \, x_k \in K, \\ &
        x \text{ and } x_k \text{ are are in different connected components of $K$ for each } k \}.
    \end{split}
\end{equation} 
Let $\vec{f} \in J^{m-1}(K)$ satisfy $f^{(j)}=0$ for $|j|\ge 1$ and $\|M^{(m)} \vec{f} \,\|_{L^p} < \infty$, and let $\eta_\eps$ be defined by \eqref{eq:def_eta_eps}. Then the family $ \vec{f}_{\varepsilon} = \{ f_\varepsilon^{(j)} \}_{|j| \le k} $ defined by
\[
f_\varepsilon^{(j)}:=
\begin{cases}
    \eta_\varepsilon \circ f^{(0)} &\text{if $j=0$}\\
    0 & \text{otherwise}
\end{cases}
\]
satisfies
\begin{equation}%\label{eq:smallness_compression_besov}
    \| M^{(m)} \vec{f}_\varepsilon \, \|_{L^p} \le \delta(\eps)
\end{equation}
for some function $\delta:\R_+\to\R_+$ with $\delta(\eps)\to 0$ as $\eps\to 0^+$.
\end{proposition}

\begin{proof}
%\giacomo{It might be useful to define the jet norm of $\vec{f}$ to be the quantity above. In this way, the main proof becomes a little more consistent}

%Since \giacomo{$m \geq 2$}, $\nabla E^{(m-1)} \vec{f} = 0$ pointwisely. Here, 
We apply \Cref{thm:extension_sobolev_p>2} to extend $\vec{f}$. Therefore, $f^{(0)}$ is constant on each connected component of $K$.

From the properties of $\eta_\varepsilon$, it is easy to see that
\begin{align}
    M^{(m)} \vec{f}_\varepsilon (x) \leq M^{(m)} \vec{f} (x),
        \quad \text{for any } \varepsilon>0 \text{ and any } x \in \R^2.
\end{align}

For any $x \in K^c$, $\dist (x,K) > 0$ and hence $M^{(m)} \vec{f}_\varepsilon (x) \rightarrow 0$ as $ \varepsilon \rightarrow 0 $.

For $x \notin S(K) \cup K^c $, we prove $M^{(m)} \vec{f}_\varepsilon (x) \rightarrow 0$ as $ \varepsilon \rightarrow 0$ by contradiction. Indeed, if there exists $c>0$ with $M^{(m)} \vec{f}_{\varepsilon_k} (x) > c$ for some sequence $\varepsilon_k \rightarrow 0$, we can find a sequence $\{(y_k,z_k)\}_k$ such that 
\begin{align}   \label{e:compression_sobolev:6}
    \frac{ | \eta_{\varepsilon_k} \circ f^{(0)}(y_k) - \eta_{\varepsilon_k} \circ f^{(0)}(z_k) | }
        { |x-y_k|^m + |x-z_k|^m } > \frac{c}{2}.
\end{align}
From the fact that $\varepsilon_k \rightarrow 0$ we deduce that $\lim_k {y_k} = \lim_k {z_k} = x$. Since $x \notin S(K) \cup K^c$ and $f^{(0)}$ is component-wise constant, we have $|\eta_{\varepsilon_k} \circ f^{(0)}(y_k) - \eta_{\varepsilon_k} \circ f^{(0)}(z_k)| = 0$ for $k$ large enough, which leads to a contradiction with \eqref{e:compression_sobolev:6}.

Note that $\|M^{(m)} \vec{f} \,\|_{L^p} < \infty$ and $|S(K)| = 0$. Combining the information above, we have that for any $x \notin S(K)$, $M^{(m)} \vec{f}_\varepsilon (x) \rightarrow 0$ as $ \varepsilon \rightarrow 0$. Using that fact $|S(K)| = 0$, by dominated convergence we conclude that $\| M^{(m)} \vec{f}_\varepsilon \, \|_{L^p} \rightarrow 0$ as $\eps\to 0$.
\end{proof}

\subsection{Trace theorem in Sobolev spaces with $p \leq n$}\label{subsec:trace_besov}

%\bian{Can you move trace and compression results for $p \leq n$ to here?}

%\newpage
%\section{\bian{This section is to be moved into other sections}}

Given a compact, $d$-regular set $K\subset \R^n$, we recall from \cite{Jonsson-Wallin} the notion of Besov space $B^{p,q}_\beta(K)$. As shown in \cite{Jonsson-Wallin}, they describe precisely the traces of Sobolev spaces onto compact, $d$-regular sets. We use $\mu$ to denote a fixed $d$-regular measure on $K$, i.e. $\mu:=\H^d|_K$. 

\begin{definition}[{\cite[Chapter~V,\S~2.3,~Definition~2]{Jonsson-Wallin}}]\label{def:Besov_space}
    Let $\beta>0$ and $k<\beta\le k+1$, with $k$ integer. The collection $\vec{f}=\{f^{(j)}\}_{|j|\le k}$ belongs to the \textit{Besov space} $B^{p,q}_\beta(K)$ if and only if there is a sequence of jets $\vec{f}_\nu=\{f_\nu^{(j)}\}_{|j|\le [\beta]}$, indexed by $\nu\in\mathbb{N}$ with $f_\nu^{(j)}\in L^p(\mu)$, and a sequence $(a_\nu)_{\nu\in\mathbb{N}} \in \ell^q(\mathbb{N}) $ such that for every $\nu\in\mathbb{N}$:
    \begin{enumerate}[label=\alph*)]
        \item $\|f^{(j)}-f^{(j)}_\nu\|_{L^p(\mu)}\le 2^{-\nu(\beta-|j|)}a_\nu$ for $|j|\le k$;
        \item $\|f^{(j)}_\nu-f^{(j)}_{\nu+1}\|_{L^p(\mu)}\le 2^{-\nu(\beta-|j|)}a_\nu=a_\nu$ if $\beta=k+1$ and $|j|=k+1$;
        \item 
        \[
        \left(2^{2\nu}\int\int_{|x-y|<2^{-\nu}} |R_{j\nu}(x,y)|^p\, d\mu(x)\,d\mu(y)\right)^{1/p}\le 2^{-\nu(\beta-|j|)}a_\nu
        \]
        for $|j|\le [\beta]$, where 
        \[
        R_{j\nu}(x,y):=f_\nu^{(j)}(x)-\sum_{|j+l|\le[\beta]}\frac{f_\nu^{(j+l)}(y)}{l!}(x-y)^l;
        \]
        \item $\|f_0^{(j)}\|_{L^p(\mu)}\le a_0$ for $|j|\le [\beta]$.
    \end{enumerate}
    The \textit{norm} of $\vec{f}$ in $B^{p,q}_\beta(K)$, denoted by $\|\vec{f}\|_{B^{p,q}_\beta(K)}$, is given by the infimum of $\left(\sum_\nu a_\nu^q\right)^{1/q}$ among all sequences $(a_\nu)_{\nu\in\mathbb{N}}$ and families $\{f_\nu^{(j)}\}_{|j|\le k}$ satisfying a)-d) above.
\end{definition}

We also recall from \cite{Jonsson-Wallin} the trace and extension theorems for Sobolev/Besov spaces. This constitutes the replacement for Fact \ref{fact:restriction} and Fact \ref{fact:extension} in the case $W^{m,p}(\R^2)$ with $p \leq 2$.

\begin{theorem}[Trace of Sobolev functions {\cite[Chapter~VII,~Theorem~1]{Jonsson-Wallin}}]\label{thm:trace_besov}
    Let $K$ be a compact $d$-regular set in $\R^n$, $0<d<n$, $1<p<\infty$, and $\beta=\alpha-(n-d)/p>0$. Then 
    \[
    W^{\alpha,p}(\R^n)|_K=B^{p,p}_\beta(K).
    \]
    More precisely:
    \begin{enumerate}[label=(\roman*)]
        \item \textbf{Restriction.}
        If $F\in W^{\alpha,p}(\R^n)$ then for every multi-index $j$ with $|j|\le \beta$, the functions $f^{(j)}:=D^jF$ are defined (and approximately continuous) $\mu$-almost everywhere, and satisfy the conditions a)-d) of Definition \ref{def:Besov_space}.
        
        \item \textbf{Extension.}
        Let $k<\beta\le k+1$, with $k$ integer. There exists an extension operator $E: B^{p,p}_\beta(K)\to W^{\alpha,p}(\R^n)$ such that, if $\vec{f}=\{f^{(j)}\}_{|j|\le k}\in B^{p,p}_\beta(K)$, then $\|Ef\|_{W^{\alpha,p}(\R^n)}\le C \|f\|_{B^{p,p}_\beta(K)}$, and $D^jF=f^{(j)}$ $\mu$-almost everywhere.
    \end{enumerate}
\end{theorem}

\begin{remark}[The case $d=n$]\label{rmk:triebel_lizorkin}
    Theorem \ref{thm:trace_besov} does not hold in general when $d=n$. However, if $p=2$ then the Sobolev space $W^{\alpha,2}$ coincides with the Lizorkin-Triebel space $F^{2,2}_\alpha$, which in turn coincides with the Besov space $B^{2,2}_\alpha$ \cite[Theoerm~4.2.2]{Adams-Hedberg}. For the latter space, a trace theorem is available also in the case $d=n$ \cite[Chapter~VI, Theorem~1]{Jonsson-Wallin}, and the trace coincides exactly with $B^{2,2}_\beta$, namely the same trace as the Sobolev space. In conclusion, by repeating the same arguments, the main result of Theorem \ref{thm:approximation_Sobolev_ple2} holds also in the case $d=2$, $p=2$.
\end{remark}

%\subsection{Compression procedure}

We also need a replacement for the compression procedure, Fact \ref{fact:norm_goes_to_0}. For this purpose, we first prove some properties of Besov functions whose derivatives vanish.

Consider $\vec{f}=\{f^{(j)}\}_{|j|\le k}\in B^{p,q}_\beta(K)$ satisfying 
\begin{equation}\label{eq:zero_derivatives_assumption}
f^{(j)}=0 \quad\text{$\mu$-almost everywhere, for every $|j|\ge 1$}.
\end{equation}
By the definition of Besov space, we know that there are approximating sequences $\{f_\nu^{(j)}\}_{|j|\le k}$, $\nu\in\mathbb{N}$, and a sequence $(a_\nu)_\nu\in \ell^q (\mathbb{N})$, with $\|\vec{f}\|_{B^{p,q}_\beta}\sim \left(\sum_\nu a_\nu^q\right)^{1/q}$. A priori, the approximating sequence might not satisfy the same property \eqref{eq:zero_derivatives_assumption}. However, the following lemma shows that we can also enforce $f_\nu^{(j)}=0$ for $|j|\ge1 $ and for every $\nu$, obtaining an equivalent norm. This result is actually formally equivalent to \cite[Chapter~V,\S~2.2, Remark~3]{Jonsson-Wallin}, but we include the proof for completeness.

%\bian{Where is the def. of being a good approx. sequence?}

\begin{lemma}[Zero derivatives]\label{lemma:zero_derivatives}
    Let $\vec{f}=\{f^{(j)}\}_{|j|\le k}\in B^{p,q}_\beta(K)$, and suppose that $f^{(j)}=0$ for every $|j|\ge 1$. Let $\{f_\nu^{(j)}\}_{|j|\le [\beta]}$ be an approximating sequence satisfying a)-d) of Definition \ref{def:Besov_space} with some constants $a_\nu$, satisfying $\sum a_\nu^q<\infty$. Then the modified sequence 
    \[
    \tilde f_\nu^{(j)}:=
    \begin{cases}
        f_\nu^{(0)} &\text{if $j=0$}\\
        0 & \text{if $|j|\ge 1$}
    \end{cases}
    \]
    satisfies a)-d) with some sequence $\tilde a_\nu$ for which $\sum_\nu\tilde a_\nu^q\le C \sum_\nu a_\nu^q$, where $C$ only depends on $\beta,p,q,\mu$.
\end{lemma}

\begin{proof}
    %By assumption $f_\nu^{(j)}$ satisfy a),b),c),d), for some sequence $a_\nu$ satisfying $\sum_\nu a_\nu^q<\infty$. We need to show that also $\tilde f_\nu^{(j)}$ satisfy the same assumptions, with a sequence $\tilde a_\nu$ comparable to $a_\nu$. 
    For $|j|\ge 1$, it is clear that all assumptions a)-d) are trivially satisfied. Indeed the remainders $\tilde R_{j\nu} $ for $\tilde f_\nu^{(j)}$, defined in \Cref{def:Besov_space}, are all zero and the other quantities at the left-hand sides of a)-d) are all zero. We thus only need to show a)-d) for $j=0$. 

    Assumption a) remains valid for $\tilde f_\nu^{(0)}$, since
    \[
    \|f^{(0)}-\tilde f_\nu^{(0)}\|_{L^p(\mu)}=\|f^{(0)}- f_\nu^{(0)}\|_{L^p(\mu)}\le 2^{-\beta \nu}a_\nu.
    \]
    
    Assumption b) is satisfied. Indeed, for $\beta<k+1$, this condition is empty. For $\beta=k+1$, if $\beta>1$, the left-hand side is zero, and if $\beta\le 1$, the left-hand side equals the one for $f_\nu$. 
    
    Assumption d) is trivially satisfied, since
    \[
    \|\tilde f_0^{(0)}\|_{L^p(\mu)}=\|f_0^{(0)}\|_{L^p(\mu)}\le a_0.
    \]
    Now we are left with checking c). The remainder simplifies to 
    \[
    \tilde R_{0\nu}(x,y)=\tilde f_\nu^{(0)}(x)-\tilde f_\nu^{(0)}(y),
    \]
    because all the higher-order terms are zero. The remainder for $\vec{f}_\nu$ is
    \[
    R_{0\nu}(x,y)=f_\nu^{(0)}(x)-f_\nu^{(0)}(y)-\sum_{1\le |l|\le \beta} \frac{f_\nu^{(l)}(y)}{l!}(x-y)^l.
    \]
    By triangle inequality
    \begin{align*}
        \bigg( 2^{d\nu}\int\int_{|x-y|<2^{-\nu}} & |\tilde R_{0\nu}(x,y)|^p\, d\mu(x)\,d\mu(y) \bigg) ^{1/p} \\ 
        &\le \left(2^{d\nu}\int\int_{|x-y|<2^{-\nu}} | R_{0\nu}(x,y)|^p\, d\mu(x)\,d\mu(y)\right)^{1/p}\\
        &+ \left(2^{d\nu}\int\int_{|x-y|<2^{-\nu}} | R_{0\nu}(x,y)-\tilde R_{0\nu}(x,y)|^p\, d\mu(x)\,d\mu(y)\right)^{1/p}.
    \end{align*}
    The first term is clearly bounded by $2^{-\nu\beta}a_\nu$ by assumption a) for $\vec{f}_\nu$. Regarding the second term, we estimate every term in the Taylor expansion,
    \begin{align}
    \sum_{1\le |l|\le \beta} &\left(2^{d\nu}\int\int_{|x-y|<2^{-\nu}} \left|\frac{f_\nu^{(l)}(y)}{l!}(x-y)^l\right|^p \, d\mu(x)\,d\mu(y)\right)^{1/p}\label{eq:double_integral}\\
    &\le \sum_{1\le |l|\le \beta}\left(2^{d\nu} \frac{1}{(l!)^p}2^{-\nu |l| p}\int\int_{|x-y|<2^{-\nu}} |f_\nu^{(l)}(y)|^p\, d\mu(x)\,d\mu(y)\right)^{1/p}\nonumber\\
    & \le\sum_{1\le |l|\le \beta}\frac{1}{(l!)}2^{-\nu |l|} \left(2^{d\nu} \int |f_\nu^{(l)}(y)|^p \mu(B(y,2^{-\nu}))\,d\mu(y)\right)^{1/p}\nonumber\\
    &\le C\sum_{1\le |l|\le \beta} 2^{-\nu |l|}\left(2^{d\nu}\int |f_\nu^{(l)}(y)|^p 2^{-d\nu}\,d\mu(y)\right)^{1/p}\nonumber\\
    & \le C\sum_{1\le |l|\le \beta}2^{-\nu |l|}\|f_\nu^{(l)}\|_{L^p(\mu)}\label{eq:last_term}
    \end{align}
    Now from b) and d) (in the case $|j|=\beta$) or from a) and d) (in the case $|j|\le k$) it follows that
    %$\|f_\nu^{(l)}\|_{L^p(\mu)}\le 2^{-\nu(\alpha-|l|)}a_\nu$ for $|l|\le [\alpha]$ . By a)/b) and d) and by telescopic sum it follows that 
    $\|f_\nu^{(l)}\|_{L^p(\mu)}\le 2\sum_{i\le\nu}  a_i$ (cf. \cite[Chapter~V,\S~2.2, Remark~2]{Jonsson-Wallin}). We deduce that the last term in \eqref{eq:last_term} is bounded by
    \[
    C\sum_{1\le |l|\le \beta}  2^{-\nu|l|} 2\sum_{i=0}^\nu a_i=: \tilde a_\nu.
    \]
    By Hardy's inequality\footnote{$\sum_{\nu=0}^\infty 2^{a\nu}\left(\sum_{i=0}^\nu a_i \right)^q\le c\sum_{\nu=0}^\infty 2^{a\nu}a_\nu^q$ for $a<0$. See \cite[Chapter V, \S2.1, Lemma~3]{Jonsson-Wallin}.}, we have
    \[
    \sum_\nu \tilde a_\nu^q \le C \sum_{1\le |l|\le \beta} \sum_\nu 2^{-\nu|l|q}\left(\sum_{i=0}^\nu a_i\right)^q\le C_2 \sum_{1\le |l|\le \beta} \sum_\nu 2^{-\nu q|l|}a_\nu^q\le C_3 \sum_\nu a_\nu^q.
    \]
    This concludes the proof.
\end{proof}

We now consider any $\vec{f}=\{f^{(j)}\}_{|j|\le k}\in B^{p,q}_\beta(K)$ satisfying $f^{(j)}=0$ for $|j|\ge 1$. We write for simplicity $f=f^{(0)}$. Using Lemma \ref{lemma:zero_derivatives}, this implies the existence of a sequence $(a_\nu)_\nu\in\ell^q(\mathbb{N})$ and an approximating sequence $f_\nu$ of $L^p(\mu)$ functions (corresponding to $f_\nu^{(0)}$) satisfying the following simplified conditions analogous to a)-d) in \Cref{def:Besov_space},
\begin{enumerate}[label=\alph*')]
    \item $\|f-f_\nu\|_{L^p(\mu)}\le 2^{-\beta\nu}a_\nu$;
    \item This condition is empty (the minimum admissible value for $|j|$ is 1, for which all functions are zero);
    \item 
    \[
    \left(2^{2\nu}\int\int_{|x-y|<2^{-\nu}} |f_\nu(x)-f_\nu(y)|^p\,d\mu(x)\,d\mu(y)\right)^{1/p}\le 2^{-2\nu}a_\nu;
    \]
    \item $\|f_0\|_{L^p_\mu(K)}\le a_0$.
\end{enumerate}

With above observation, we are ready to prove the compression property for Besov spaces on $K$.

\begin{proposition}[Monotone compression in Besov spaces]\label{prop:compression_besov}
Consider the map $\eta_\eps$ from \eqref{eq:def_eta_eps}. Let $\vec{f}=\{f^{(j)}\}_{|j|\le k}\in B^{p,q}_\beta(K)$ satisfy $f^{(j)}=0$ for $|j|\ge 1$. Then the family $\vec{f}_\eps=\{ f_\eps^{(j)}\}_{|j|\le k}$ defined by
\[
 f_\eps^{(j)}:=
\begin{cases}
    \eta_\eps\circ f^{(0)} &\text{if $j=0$}\\
    0 & \text{otherwise}
\end{cases}
\]
also belongs to $B^{p,q}_\beta(K)$. Moreover, it satisfies
\begin{equation}\label{eq:smallness_compression_besov}
\|\vec{f}_\eps\|_{B^{p,q}_\beta(K)}\le \delta(\eps)
\end{equation}
for some function $\delta:\R_+\to\R_+$ with $\delta(\eps)\to 0$ as $\eps\to 0^+$.
\end{proposition}

\begin{proof}
    As $\vec{f}$ belongs to the Besov space $B^{p,q}_\beta$, we can find sequences $(f_\nu)_\nu$ and $(a_\nu)_\nu$ that satisfy conditions a)-d) of Definition \ref{def:Besov_space}. We claim that the sequence $ f_{\eps,\nu}:=\eta_\eps\circ f_{\nu}$ is a valid approximating sequence for $f_\eps$, i.e., satisfying all the conditions specified in \Cref{def:Besov_space}. Indeed observe that, by the compression property, $| f_{\eps,\nu}^{(j)}- f_{\eps,\nu+1}^{(j)}|\le |f_\nu^{(j)}-f_{\nu+1}^{(j)}|$ and $| f_\eps^{(j)}-f_{\eps,\nu}^{(j)}|\le |f^{(j)}-f_\nu^{(j)}|$ pointwise. It follows that $\vec{f}_\eps\in B^{p,q}_\beta(K)$, since the left-hand sides in a)-d) do not increase after composing with $\eta_\eps$. In particular this shows that $(a_\nu)_\nu$ is a valid sequence also for the estimates on $ f_\eps$ specified in \Cref{def:Besov_space}.
    
    Next we use the smallness assumption on the $C^0$ norm of $\eta_\eps$ to show \eqref{eq:smallness_compression_besov}, by estimating each term in a)-d) of \Cref{def:Besov_space}.
    
    Regarding the term in a), we have that $\| f_\eps-f_{\eps,\nu}\|_{L^p_\mu(K)}\le 2\eps \mu(K)^{1/p}$, due to $\| f_\eps\|_\infty \le \eps$ and $\| f_{\eps,\nu}\|_\infty \le \eps$. By the compression property, $\| f_\eps-f_{\eps,\nu}\|_{L^p_\mu(K)}\le \| f-f_\nu\|_{L^p_\mu(K)} \le 2^{-2\nu}a_\nu$. It follows that
    \[
    2^{2\nu}\| f_\eps- f_{\eps,\nu}\|_{L^p_\mu(K)}\le \min\{ 2^{2\nu+1}\eps \mu(K)^{1/p},a_\nu\}.
    \]

    The condition b) is empty, since the minimum admissible value for $|j|$ is 1, for which all $f_{\eps,\nu}^{(j)}$ are zero.
    
    The term in d) is similar to a) above, i.e. $\| f_{\eps,0}\|_{L^p_\mu(K)}\le \min\{ \eps\mu(K)^{1/p},a_0\}$. 
    
    Regarding the term in c), we have two estimates: on the one hand by the compression property
    \begin{align}   \label{e:compression_besov:8}
    \left(2^{2\nu}\int\int_{|x-y|<2^{-\nu}} | f_{\eps,\nu}(x)- f_{\eps,\nu}(y)|^p\,d\mu(x)\,d\mu(y)\right)^{1/p}\le 2^{-2\nu}a_\nu.
    \end{align}
    On the other hand, since $\| f_{\eps,\nu}\|_\infty\le \eps$, we have
    \begin{align*}
    &\left(2^{2\nu}\int\int_{|x-y|<2^{-\nu}} | f_{\eps,\nu}(x)- f_{\eps,\nu}(y)|^p\,d\mu(x)\,d\mu(y)\right)^{1/p} \\
    &\le \left(2^{2\nu}(2\eps)^p\int\,d\mu(x)\int_{|x-y|<2^{-\nu}} \,d\mu(y)\right)^{1/p}\\
    &\le C\eps \left(2^{2\nu}\mu(K)2^{-2\nu}\right)^{1/p}\\
    &\le C\eps.
    \end{align*}
    Here we used the upper regularity of the measure $\mu$ to deduce that $\mu(B(y,2^{-\nu}))\le C 2^{-d\nu}$.
    In conclusion, the left-hand side of \eqref{e:compression_besov:8} is bounded by $\min\{2^{-2\nu}a_\nu,C\eps\}$.

    With these observations, we define the sequence 
    \begin{align*}
     a_{\eps,0}&:=\min\{a_0,\eps \mu(K)^{1/p}\}, \\ 
     a_{\eps,\nu}&:=\max\left\{\min\{a_\nu, C 2^{2\nu}\eps\},\min\{ 2^{2\nu+1}\eps \mu(K)^{1/p},a_\nu\}\right\},\qquad \nu\ge 1.
    \end{align*}
    Then we have
    \begin{enumerate}[label=\alph*')]
    \item $\| f_\eps- f_{\eps,\nu}\|_{L^p_\mu(K)}\le 2^{-2\nu}  a_{\eps,\nu}$
    \item $\| f_{\eps,0}\|_{L^p_\mu(K)}\le a_{\eps,0}$ \setcounter{enumi}{3}
    %\item This condition is automatically satisfied.
    \item
    \[
    \quad \left(2^{2\nu}\int\int_{|x-y|<2^{-\nu}} | f_{\eps,\nu}(x)-\tilde f_{\eps,\nu}(y)|^p\,d\mu(x)\,d\mu(y)\right)^{1/p}\le 2^{-2\nu} a_{\eps,\nu}.
    \]
    \end{enumerate}
    This shows that $ a_{\eps,\nu}$ is a valid sequence to bound the estimates for $ f_\eps$ in \Cref{def:Besov_space}, as it satisfies conditions a)-d). It remains to see that, as $\eps\to 0$, the dominated convergence theorem implies that 
    \[
    \left(\sum_\nu  a_{\eps,\nu}^q\right)^{1/q}\to 0.
    \]
    Since the last sum estimates from above $\|\vec{f}_\eps\|_{B^{p,q}_\beta(K)}$ up to some uniform constant, this concludes the proof.
\end{proof}

\section{Proof of the main approximation theorems} \label{ms:main_proof}

In this section we prove the main results for Sobolev spaces: Theorem \ref{thm:approximation_Sobolev_p>2}, Theorem \ref{thm:approximation_Sobolev_ple2} and Corollary \ref{cor:approximation_Sobolev_ple2}.
%and Theorem \ref{thm:approximation_Cmgamma}.
The strategy is the same as the one adopted in \Cref{ms:quicktour} for the $C^1$ case. However, we now use the results proven in the previous section as replacements for Facts 1-4. We will treat the cases $p>2$ and $p\le 2$ in parallel. Case 1 refers to $p>2$, Case 2 refers to $p\le 2$, with Subcase 2.I being $2-p<d\le 2$ and Subcase 2.II being $d\le 2-p$.

Moreover, as already anticipated, we will give the full proof of the case $m=1$ only (corresponding to the initial vector field $u$ belonging to $W^{1,p}$). The proof for general $m$ is virtually identical, and we comment on this in Remark \ref{remark:Wmp_difference}.

\begin{proof}[Proof of \Cref{thm:approximation_Sobolev_p>2} and \Cref{thm:approximation_Sobolev_ple2} for $m=1$]

We follow the steps from Section \ref{ms:quicktour}, commenting on the required modifications. Observe first that if $d=0$ then the set $\Omega^c$ is discrete and has a finite number of connected components. Thus Theorem \ref{thm:approximation_Sobolev_ple2} follows from Theorem \ref{thm:sverak}. Hence, from here on, we can assume that $d>0$ in Case 2 ($p\le 2$).

\begin{enumerate}[label=(S\arabic*)]
    \item \textit{Potential.} Since $\divr u=0$ and $u\in W^{1,p}(\R^2;\R^2)$, we can consider a potential $F$ for $u$, namely, $F\in W^{2,p}_\loc(\R^2)$ such that $\nabla^\perp F=u$. By the Morrey-Sobolev embedding theorem, $F$ admits a continuous representative, and we assume that $F$ coincides with this representative. Therefore, $F$ is well-defined at every point. Moreover, since $\nabla F=u^\perp=0$ on $\Omega^c$, without loss of generality, we can assume that $F$ is compactly supported.

    \textit{Case 1} ($p>2$).
    By the Morrey-Sobolev embedding $\nabla F$ is continuous, and we have that $\nabla F=0$ pointwisely on $K$. 

    \textit{Case 2} ($p\le 2$).
    In this case, $\nabla F=0$ $C_{1,p}$-quasi everywhere on $K$. 
    
    \textit{Subcase 2.I} ($2-p<d\le 2$). By Theorem \ref{thm:capacity_vs_hausdorff} we have that $\nabla F=0$ $\mu$-a.e. on $K$.
    
    \textit{Subcase 2.II} ($d\le 2-p$). In this case we have no restriction on $\nabla F$.

    \item \textit{Restriction.} 

        \textit{Case 1} ($p>2$). We apply \Cref{thm:restriction_sobolev_p>2} to obtain the restriction $\vec{f}=\{f^{(j)}\}_{|j|\le 1}$ to the set $K$. Observe that $F$ and $\nabla F$ are continuous by the Morrey-Sobolev embedding. In particular, all the $f^{(j)}$ are continuous.

        \textit{Case 2} ($p\le 2$).
        We apply the restriction part in \Cref{thm:trace_besov} with $\alpha=2$. %Observe that $F$ is continuous, since it belongs to $W^{2,p}(\R^2)$.

        \textit{Subcase 2.I} ($2-p<d\le 2$). In this case $\beta=2-\frac{2-d}{p}\in(1,2]$, and the restriction $\vec{f}=\{f^{(j)}\}_{|j|\le 1}:=\{(D^jF)|_K\}_{|j|\le 1}$ belongs to $B^{p,p}_\beta(K)$. Moreover $f^{(j)}=0$ $C_{1,p}$-quasi everywhere for $|j|=1$, and thus also $\mu$-almost everywhere by Theorem \ref{thm:capacity_vs_hausdorff}. According to Definition \ref{def:Besov_space} with $k=1$, we have that $\vec{f}=\{f^{(j)}\}_{|j|\le 1}$ satisfies a)-d) of \Cref{def:Besov_space}, namely there exist sequences $(f^{(j)}_\nu)_\nu\subset L^p(\mu)$, $|j|\le 1$, and $(a_\nu)_\nu\subset \R$ such that
        \begin{enumerate}
            \item $\|f^{(j)}-f^{(j)}_\nu\|_{L^p(\mu)}\le 2^{-\nu\beta}a_\nu$;
            \item $\|f^{(j)}_\nu-f^{(j)}_{\nu+1}\|_{L^p(\mu)}\le 2^{-\nu\beta}a_\nu$ if $\beta=2$ and $|j|=2$;
            \item 
            \[
            \left(2^{2\nu}\int\int_{|x-y|<2^{-\nu}} |f^{(0)}_\nu(x)-f^{(0)}_\nu(y)|^p\, d\mu(x)\,d\mu(y)\right)^{1/p}\le 2^{-\nu\beta}.
            \]
            Here due to \Cref{lemma:zero_derivatives}, we also choose without loss of generality $f^{(j)}_\nu=0$ for $|j|=1$.
            \item $\|f_0^{(j)}\|_{L^p(\mu)}\le a_0$ for $|j|\le[\beta]$.
        \end{enumerate}

        \textit{Subcase 2.II} ($d\le 2-p$). In this case $\beta=2-\frac{2-d}{p}\in (0,1]$. Therefore, due to \Cref{thm:trace_besov}, the restriction jet is given by $\vec{f}=\{f^{(0)}\}:=\{F|_K\}$, satisfying the same conditions above.
        
    \item \textit{Sard.} %\giacomo{$K$ itself is not compact. Since $K = \Omega^c$ for some bounded domain $\Omega$, $K$ is the union of an unbounded closed set and a compact set. }
        By Lemma \ref{lemma:sard_powered}, we have that $F(K)$ is a compact set with zero measure. As in the $C^1$ case, for every $\eps>0$, we can therefore find finitely many intervals $\{I_i^\eps\}_{1 \leq i \leq N}$ of total measure at most $\eps$ and the well-separated open sets $\{U_i^\eps\}_{1 \leq i \leq N}$ whose union contains $\Omega^c$. Here, note that $F$ is constant on $\tilde K$ in Case 2.

    \item \textit{Monotone compression.}
        There is no essential modification compared to the $C^1$ case. We define the compression map $\vec{f}_\eps:=\eta_\eps\circ \vec{f}$, where $\eta_\eps$ is defined as in \eqref{e:eta_eps_qt}. This means that $f_\eps^{(0)}=\eta_\eps\circ f^{(0)}$ and $f_\eps^{(j)}=0$ for $|j|=1$.
        
    \item \textit{Estimate on the compressed jets.} \label{it:compression_estimate}
    
        \textit{Case 1} ($p>2$). With $|S(\Omega^c)| = 0$ and Proposition \ref{prop:compression_sobolev}, we deduce that the compressed jets $\vec{f}_\eps$ satisfy $\|M^{(m)}\vec{f}_\eps\|_{L^p}\le \delta(\eps)$, where $\delta(\eps)\to 0$ as $\eps \to 0^+$.

        \textit{Case 2} ($p\le 2$). We apply Proposition \ref{prop:compression_besov} to deduce that the compressed jets $\vec{f}_\eps$ satisfies $\|\vec{f}_\eps\|_{B^{p,p}_\beta(K)}\le \delta(\eps)$, where $\delta(\eps)\to 0$ as $\eps\to 0^+$. This is true both in Subcase 2.I and in Subcase 2.II.
        
    \item \textit{Extension.}

    \textit{Case 1} ($p>2$).
        We first apply \Cref{thm:extension_sobolev_p>2} with the compact set given by $\Omega^c \cap \bar{B}_R(0)$ with sufficiently large $R>0$ such that $\bar \Omega \subset B_R(0)$, thus we extend $\vec{f}_\eps$ to a function $F_\eps\in W^{2,p}(\R^2)$ with $\nabla F_\eps = 0$ on $\Omega^c$. From the definition of the Whitney extension operator, $F_\eps$ is constant on each connected component of $\Omega^c$. By the estimates in \ref{it:compression_estimate} and in \Cref{thm:extension_sobolev_p>2}, we have $\|\nabla^2 F_\eps\|_{L^p(\R^2)} \lesssim \delta(\eps)$. 
        
        %Then we multiply by the cutoff $\theta$. By the Poincare' inequality and standard estimates it also follows that $\|\theta F_\eps\|_{W^{2,p}(\R^2)}=o(\eps)$. 
        
        %Strictly speaking, in Theorem \ref{thm:extension_sobolev_p>2} only $\|\nabla^2 F\|_{L^2}$ is controlled, but since $K$ is compact we can also assume that $F$ is compactly supported (up to multiplying by a universal cutoff), and then by the Poincar\'{e} inequality we can control the full Sobolev norm $\|F\|_{W^{2,p}(\R^2)}$.

    \textit{Case 2} ($p\le 2$).   
        Recall that the bounded domain $\Omega$ admits the decomposition $\Omega^c=K\cup \tilde K$, with $\tilde K$ being the unbounded component and $K$ being $d$-regular. We fix a cutoff function $\chi \in C^\infty_c(\R^2)$ such that $\chi$ is supported on $\tilde K ^c$ and such that $\chi \equiv 1$ on $K$.
        
        We apply the extension part in Theorem \ref{thm:trace_besov} to extend $\vec{f}_\eps$ to a function $F^*_\eps \in W^{2,p}(\R^2)$. Define $F_\eps:= \chi F^*_\eps$. From the estimates in \ref{it:compression_estimate} we deduce that $\|F_\eps\|_{W^{2,p}(\R^2)} \lesssim C(\chi) \delta(\eps)$. Here, $C(\chi)>0$ is a constant depending on $\chi$ and independent of $\eps$.
        
        In Subcase 2.I, $F_\eps$ satisfies $D^j F_\eps=f_\eps^{(j)}$ $\mu$-a.e. on $K$, for every $|j|\le 1$. However, since $F_\eps$ belongs to $W^{2,p}(\R^2)$ it is also continuous, hence we conclude that $F_\eps=f_\eps^{(0)}$ everywhere on $K$.

        In Subcase 2.II, the function satisfies $F_\eps=f_\eps^{(0)}$ $\mu$-a.e. on $K$. However, for the same reason above $F_\eps=f_\eps^{(0)}$ everywhere on $K$.
        
    \item \textit{Auxiliary function.}

        In this step we construct an auxiliary function $h_\eps \in C^\infty(\R^2)$ with
        \begin{align}
            \supp \nabla h_\eps \subset&\, \Omega,    \label{e:main_proof:8} \\ 
            g_\eps := F - F_\eps-h_\eps =&\, 0, \quad \text{on } \Omega^c.
                   \label{e:main_proof:10}
        \end{align}
    
        \textit{Case 1} ($p>2$).
        Due to the definitions of $\eta_\eps$, $U_i^\eps$, and the same reason as the $C^1$ case in \Cref{ms:quicktour}, $F-F_\eps = c_i^\eps$ on $\Omega^c \cap U_i^\eps$ for any $i$. $\{ \Omega^c \cap U_i^\eps \}_{1 \leq i \leq N}$ are finitely many disjoint closed sets, hence we can find $h_\eps\in C^\infty$ such that $h_\eps = c_i^\eps$ on $U_i^\eps$ and such that \eqref{e:main_proof:8} and \eqref{e:main_proof:10} hold.

        Furthermore, we have $\nabla g_\eps = 0$ on some neighborhood of $\Omega^c$.

        \textit{Case 2} ($p\le 2$). 
        For the same reason above, $F-F_\eps = c_i^\eps$ on $K\cap U_i^\eps$ for any $i$. Moreover, $F-F_\eps = \tilde c$ on $\tilde K$. Similarly, we can find $h_\eps\in C^\infty$ such that $h_\eps = c_i^\eps$ on $U_i^\eps$, $h_\eps = \tilde c$ on $\tilde K$ and \eqref{e:main_proof:8}, \eqref{e:main_proof:10} hold.

        In Subcase 2.I, we also have $\nabla g_\eps = 0$ for $\mu$-a.e. on $K$, and thus by Theorem \ref{thm:capacity_vs_hausdorff} also $C_{1,p}$-quasi everywhere on $K$. 
    
    \item \textit{Hedberg's theorem and conclusion.} 
        Now we apply Theorem \ref{thm:hedberg} to find a sequence $g_\eps^k\in C^\infty_c( \Omega )$ with $\|g_\eps^k-g_\eps\|_{W^{2,p}}\to 0$ as $k\to\infty$, and we define the sequence $u_\eps^k:=\nabla^\perp(g_\eps^n+h_\eps)$. Since $\nabla^\perp(g_\eps^k+h_\eps)=0$ on some neighborhood of $\Omega^c$, we have $u_\eps^k\in C^\infty_c(\Omega)$. Also, we can estimate
    \begin{align*}
        \|u_\eps^k-u\|_{W^{1,p}}&=\|\nabla^\perp(g_\eps^k+h_\eps-F)\|_{W^{1,p}}\\
        &\le \|\nabla^\perp(g_\eps+h_\eps-F)\|_{W^{1,p}}+\|\nabla^\perp(g_\eps-g_\eps^k)\|_{W^{1,p}}\\
        &=\|\nabla^\perp F_\eps\|_{W^{1,p}}+\|\nabla^\perp(g_\eps^k-g_\eps)\|_{W^{1,p}}\\
        &\le \|F_\eps\|_{W^{2,p}}+\|g_\eps^k-g_\eps\|_{W^{2,p}}.
    \end{align*}
        Now the first summand goes to zero as $\eps\to 0$ by Step (S6), while the second summand goes to zero as $k\to\infty$ by the conclusion of \Cref{thm:hedberg}. Taking a diagonal sequence, we can extract an approximating sequence of divergence-free vector fields $\{u_k\}_k$ converging to $u$, which concludes the proof.\qedhere
\end{enumerate}
\end{proof}

\begin{remark}[\Cref{thm:approximation_Sobolev_p>2} and \Cref{thm:approximation_Sobolev_ple2} for $m\geq 2$]    \label{remark:Wmp_difference}
%\giacomo{I added this remark about Sobolev spaces of arbitrary order. I think it is sufficient, because the same proof applies as explained below.}

If in Theorem \ref{thm:approximation_Sobolev_p>2} and Theorem \ref{thm:approximation_Sobolev_ple2} $u$ is assumed to be of class $W^{m,p}$, then the potential $F$ belongs to $W^{m+1,p}_\loc(\R^2)$, and virtually the same proof applies. The only changes are in the order of the derivatives that we are considering in each step.
In particular, by the Morrey-Sobolev embedding all derivatives of $F$ up to order $m-1$ are continuous. The role that $\nabla F$ has in the proof above is now taken by $\nabla^m F$.

In Case 1 ($p>2$), $\nabla^{m}F$ is continuous. 
In Case 2 ($p\le 2$), we have $\nabla^m F=0$ $C_{1,p}$-quasi everywhere. As a consequence we have that $\nabla^m F=0$ $\mu$-a.e. in Subcase 2.I ($2-p<d\le 2$), while in Subcase 2.II ($d\le 2-p$) we have no restriction on $\nabla^m F$.

Regarding the trace theorems (Theorems \ref{thm:extension_sobolev_p>2} and \ref{thm:restriction_sobolev_p>2} for $p>2$, and Theorem \ref{thm:trace_besov} for $p\le 2$), they are stated for Sobolev spaces of arbitrary integer order and we can apply them directly. The trace space becomes the Besov space $B^{p,p}_\beta$ with $\beta=m-\frac{2-d}{p}\in (m-1,m+1]$. Similarly, the monotone compression procedure of Proposition \ref{prop:compression_besov} works for arbitrary order.

Finally, Hedberg's result also holds for Sobolev spaces of any order, concluding the proof for $m\ge2$.
\end{remark}

We now turn to the proof of of Corollary \ref{cor:approximation_Sobolev_ple2}. This will be a consequence of the more general lemma below. We say that an open set $\Omega\subset\R^n$ has the Sobolev approximation property if $W^{m,p}_{0,\divr}(\Omega)=\tilde W^{m,p}_{0,\divr}(\Omega)$ (recall \eqref{e:intro:6}-\eqref{e:intro:8}).

%the conclusion corresponding to Theorem \ref{thm:approximation_Sobolev_p>2} or Theorem \ref{thm:approximation_Sobolev_ple2} holds for every divergence-free vector field $u\in W^{m,p}$ satisfying the assumptions of the theorems. 

\begin{lemma}\label{lemma:intersection_property}
    Let $\Omega\subset \R^2$ be an open bounded set, such that $\Omega^c=\tilde K\cup K_1\cup K_2$, with $K_1$, $K_2$ and $\tilde K$ pairwise disjoint, $K_1$ and $K_2$ compact sets, and $\tilde K$ being the only unbounded connected component of $\Omega^c$. Suppose that the Sobolev approximation property holds for $\Omega_1:=\Omega\cup K_1$ and $\Omega_2:=\Omega\cup K_2$. Then it holds for $\Omega$.
\end{lemma}

\begin{proof}
    Let us fix a cutoff function $\chi\in C^\infty_c(\Omega)$ such that $\chi\equiv 1$ on a neighborhood of $K_1$ and $\supp \chi\subset K_2^c$. Observe that $1-\chi\equiv 1$ on a neighborhood of $K_2$ and $\supp \,(1-\chi)\subset K_1^c$.
    
    Consider now a divergence-free vector field $u\in W^{m,p}(\R^2;\R^2)$ with $D^j u=0$ on $\Omega^c$ for $1\le |j|\le m-1$, and fix $\eps>0$. Since the approximation property holds for $\Omega_1$ and $\Omega_2$ separately, we can find divergence-free smooth vector fields $u_1$ and $u_2$ with
    \[
    \supp u_i\subset (\Omega\cup K_i)^c,\qquad \|u-u_i\|_{W^{m,p}}\le \eps\quad\text{for $i=1,2$.}
    \]
    Let us call $\phi_1$, $\phi_2$ and $\phi$ the compactly-supported potentials relative to $u_1$, $u_2$ and $u$ respectively, namely
    \[
    \nabla^\perp\phi=u,\qquad \nabla^\perp \phi_i=u_i\quad \text{for $i=1,2$.}
    \]
    Observe that by the Poincar\'{e} inequality we also have
    \[
    \|\phi-\phi_i\|_{W^{m+1,p}}\le C(\Omega)\eps\quad\text{for $i=1,2$.}
    \]
    Defining now the potential $\psi:=\chi \phi_1+(1-\chi)\phi_2$, we claim that the vector field $v:=\nabla^\perp \psi$ approximates $u$ up to order $\eps$ and has the desired properties. Indeed we can write
    \[
    \phi-\psi=(\phi-\phi_1)\chi+(\phi-\phi_2)(1-\chi),
    \]
    from which it follows that
    \[
    \|\phi-\psi\|_{W^{m+1,p}}\le C(\Omega,\chi)\big(\|\phi-\phi_1\|_{W^{m+1,p}}+\|\phi-\phi_2\|_{W^{m+1,p}}\big).
    \]
    In particular,
    \[
    \|u-v\|_{W^{m,p}}\le  C_1(\Omega,\chi) \eps.
    \]
    Moreover $v$ is divergence-free by construction. It remains to check that $\psi$ is constant on a neighborhood of $\Omega^c$, and thus $v$ is supported on $\Omega^c$. To see this, observe that on a sufficiently small neighborhood of $K_1$ the function $(1-\chi)$ is zero, the potential $\phi-\phi_1$ is constant and $\chi$ attains value 1, so that $\psi$ is constant there. An analogous thing happens on a neighborhood of $K_2$. Moreover $\supp(\psi)\subset \tilde K$ by construction. This concludes the verification of all the properties and thus the Sobolev approximation property holds for $\Omega$.
\end{proof}

\begin{proof}[Proof of Corollary \ref{cor:approximation_Sobolev_ple2}]
    By either Theorem \ref{thm:approximation_Sobolev_ple2} or Theorem \ref{thm:sverak}, we know that the Sobolev approximation property holds when $\Omega^c=\tilde K\cup K$ with $K$ either connected or $d$-Ahlfors regular for some $d\in [0,2)$.  It remains to apply repeatedly Lemma \ref{lemma:intersection_property} to reach the conclusion.
\end{proof}

\bibliographystyle{plain}
\bibliography{biblio}

\end{document}